\documentclass[11pt,english]{amsart}
\usepackage[latin1]{inputenc}
\usepackage{geometry}
\usepackage{setspace}
\usepackage{amssymb}

\makeatletter

\providecommand{\LyX}{L\kern-.1667em\lower.25em\hbox{Y}\kern-.125emX\@}

 \theoremstyle{plain}    
 \newtheorem{thm}{Theorem}[section]
 \numberwithin{equation}{section} %% Comment out for sequentially-numbered
 \numberwithin{figure}{section} %% Comment out for sequentially-numbered
 \theoremstyle{remark}    
 \newtheorem*{acknowledgement*}{Acknowledgement} 
 \theoremstyle{plain}    
 \newtheorem{cor}[thm]{Corollary} %%Delete [thm] to re-start numbering
 \theoremstyle{definition}
 \newtheorem{defn}[thm]{Definition}
 \theoremstyle{definition}
 \newtheorem*{defn*}{Definition}
 \theoremstyle{plain}    
 \newtheorem{lem}[thm]{Lemma} %%Delete [thm] to re-start numbering
 \theoremstyle{plain}    
 \newtheorem{prop}[thm]{Proposition} %%Delete [thm] to re-start numbering
 \theoremstyle{remark}
 \newtheorem*{rem*}{Remark}
 \theoremstyle{plain}    
 \newtheorem*{thm*}{Theorem} 

\usepackage{setspace}

\theoremstyle{plain}
\newtheorem*{question}{Question}
\newtheorem*{proposition*}{Proposition}

\DeclareMathOperator{\dom}{dom}

\usepackage{babel}
\makeatother

\makeatother

\usepackage{babel}
\makeatother
\begin{document}

\title{Return times, recurrence densities and entropy for actions of some
discrete amenable groups}

\author{Michael Hochman}

\email{mhochman@math.huji.ac.il}

\curraddr{Einstein Institute of Mathematics, Edmond J. Safra Campus, Givat
Ram, The Hebrew University of Jerusalem Jerusalem, 91904, Israel }

\begin{abstract}
It is a theorem of Wyner \& Ziv and Ornstein \& Weiss that if one
observes the initial $k$ symbols $X_{0},\ldots ,X_{k-1}$ of a typical
realization of a finite valued ergodic process with entropy $h$,
the waiting time until this sequence appears again in the same realization
grows asymptotically like $2^{hk}$ \cite{OW93,WyZi89}. A similar
result for random fields was obtained in \cite{OW02}: in this case
one observes cubes in $\mathbb{Z}^{d}$ instead of initial segments. 

In the present paper we describe generalizations of this. We examine
what happens when the set of possible return times is restricted:
Fix an increasing sequence of sets of possible times $\{W_{n}\}$
and define $R_{k}$ to be the first $n$ such that $X_{0},\ldots ,X_{k-1}$
recurs at some time in $W_{n}$. It turns out that $|W_{R_{k}}|$
cannot drop below $2^{hk}$ asymptotically. We obtain conditions on
the sequence $\{W_{n}\}$ which ensure that $|W_{R_{k}}|$ is asymptotically
equal to $2^{hk}$.

We consider also recurrence densities of initial blocks and derive
a uniform Shannon-McMillan-Breiman theorem: Informally, if $U_{k,n}$
is the density of recurrences of the block $X_{0},\ldots ,X_{k-1}$
in $X_{-n},\ldots ,X_{n}$ then $U_{k,n}$ grows at a rate of $2^{hk}$,
uniformly in $n$. We examine the conditions under which this is true
when the recurrence times are again restricted to some sequence of
sets $\{W_{n}\}$

The above questions are examined in the general context of finite-valued
processes parameterized by discrete amenable groups. We show that
many classes of groups have time-sequences $\{W_{n}\}$ along which
return times and recurrence densities behave as expected. An interesting
feature here is that this can happen also when the time sequence lies
in a small subgroup of the parameter group. 
\end{abstract}
\maketitle

\section{\label{sec:intro} Introduction }

\pagestyle{myheadings}

\markboth{Michael Hochman}{Return times, recurrence densities and entropy}

\subsection{Background}

Let $\{X_{n}\}_{n\in \mathbb{Z}}$ be an ergodic process with values
in some finite set $\Sigma $, defined on a probability space $(\Omega ,\mathcal{F},P)$.
Let $h$ denote the entropy of the process. For a realization $x=(x_{n})_{n\in \mathbb{Z}}$
of the process and for each $k$ consider the finite word $x_{0}x_{1}\ldots x_{k-1}$
obtained by viewing $x$ through the {}``window'' consisting of
coordinates $0,1,\ldots ,k$. The first backward recurrence time of
this word is defined as\[
R_{k}^{-}(x)=\min \{n\, :\, n>0\textrm{ and }x_{-n+i}=x_{i}\textrm{ for all }0\leq i<k\}\]
 Set $R_{k}^{-}(x)=\infty $ if the sequence $x_{0}x_{1}\ldots x_{k-1}$
does not appear again in the past portion of $x$.

In \cite{WyZi89}, A. Wyner and J. Ziv proved that $\frac{1}{k}\log R_{k}^{-}\rightarrow h$
in probability and $\limsup _{k}\frac{1}{2k+1}\log R_{k}^{-}\leq h$
almost surely. This result was strengthened by D. Ornstein and B.
Weiss, who proved in \cite{OW93} that in fact\begin{equation}
\lim _{k\rightarrow \infty }\frac{1}{k}\log R_{k}^{-}(\omega )=h\label{eq:Z-recurrence-time-theorem}\end{equation}
 almost surely. 

Wyner and Ziv's definition of $R_{k}^{-}$ uses the past and future
in an asymmetric fashion, but in \cite{OW93} Ornstein and Weiss also
define a more symmetric quantity, and prove a related limit result,
which they later generalized in \cite{OW02} to the case of $\mathbb{Z}^{d}$
processes, or so-called random fields. These are families of random
variables $\{X_{u}\}_{u\in \mathbb{Z}^{d}}$ which are stationary
with respect to {}``time shifts'' in $\mathbb{Z}^{d}$. In this
setting, the patterns we are trying to match, and the sets in which
we look for recurrences of patterns, are based on the cubes $F_{n}=[-n;n]^{d}$
(here and throughout, we denote the integer segment $\left\{ i,i+1\ldots ,j\right\} $
for $i<j$ by $[i;j]$). For a realization $x=(x_{u})_{u\in \mathbb{Z}^{d}}$
of the process, we obtain a {}``pattern'' on the cube $F_{k}$,
which we denote $x(F_{k})$, by coloring each $v\in F_{k}$ with the
color $x_{v}$. The first recurrence time of the $F_{k}$-pattern
$x(F_{k})$ is defined to be the first index $n$ for which we observe
the pattern $x(F_{k})$ centered at some point in $F_{n}$ other than
the origin. To be precise, \[
R_{k}(\omega )=\min \left\{ n\, \left|\begin{array}{c}
 \textrm{for some }v\neq 0\textrm{ in }[-n;n]^{d}\\
 x_{v+u}=x_{u}\textrm{ for every }u\in F_{k}\end{array}
\right.\right\} \]
 or $R_{k}(\omega )=\infty $ if the pattern does not recur. Then
according to \cite{OW02}, with probability $1$ \begin{equation}
\lim _{k\rightarrow \infty }\frac{d}{(2k+1)^{d}}\log R_{k}(\omega )=h\label{eq:return-time-for-Z-infinity}\end{equation}
(in \cite{OW02} there is also a proof for patterns based on the cubes
$[0;n]^{d}$).

Our goal in this paper is to describe several generalizations of these
phenomena for processes indexed by other amenable groups, though many
of the results are new even for $\mathbb{Z}$-processes

\subsection{Notation and preliminaries}

Recall that a discrete countable group $G$ is amenable if there exists
a sequence of finite subsets $\left\{ F_{n}\right\} $ of $G$ which
are asymptotically invariant, ie for all $g\in G$,\[
\lim _{n\rightarrow \infty }\frac{|F_{n}\, \Delta \, gF_{n}|}{|F_{n}|}=0.\]
 For example, the sequence $\left\{ [-n;n]^{d}\right\} _{n=1}^{\infty }$
in $\mathbb{Z}^{d}$ has this property. Such sequences are called
F\o{}lner sequences. Asymptotic invariance is equivalent to the following
statement: for every finite $K\subseteq G$, it holds that \[
\lim _{n\rightarrow \infty }\frac{1}{|F_{n}|}\left|\left\{ f\in F_{n}\, :\, Kf\subseteq F_{n}\right\} \right|=1.\]
 There are other equivalent definitions of amenability, some of which
have a more analytic flavor, but these will not interest us here;
see eg \cite{Green69}.

For $G$ a discrete amenable group we consider an ergodic process
$\{X_{g}\}_{g\in G}$ taking values in a finite set $\Sigma $, defined
on a standard probability space $(\Omega ,\mathcal{F},P)$. We will
always assume that the process arises from an ergodic measure preserving
free left action of $G$ on $\Omega $, denoted $(g,\omega )\mapsto g\omega $,
and from a measurable function $X:\Omega \rightarrow \Sigma $ such
that $X_{g}(\omega )=X(g\omega )$. We denote by $h$ the entropy
of the process.

It is convenient to associate the elements $\omega \in \Omega $ with
the $\Sigma $-colorings of $G$ they induce: $g\mapsto X_{g}(\omega )$.
We will often denote elements of $\Omega $ by the letter $x$ to
emphasize this point of view, and write $x(g)$ for $X_{g}(x)$. For
$x\in \Omega $ and $E\subseteq G$ we denote by $x(E)$ the coloring
of $E$ by colors from $\Sigma $ given by $g\mapsto x(g)$, and let
$[x(E)]\in \mathcal{F}$ be the atom \[
[x(E)]=\{x'\in \Omega \, :\, \forall f\in E\; \; x'(f)=x(f)\}.\]
 More generally, if $\varphi \in \Sigma ^{E}$ is a coloring of $E$
by $\Sigma $, we write $[\varphi ]$ for the atom defined by $\varphi $
\[
[\varphi ]=\left\{ x\in \Omega \, :\, \forall f\in E\, \; x(f)=\varphi (f)\right\} \]
 and for $\Phi \subseteq \Sigma ^{F_{n}}$ we set \[
[\Phi ]=\bigcup _{\varphi \in \Phi }[\varphi ].\]

It is known that the ergodic and entropy theorems hold along certain
F\o{}lner sequences. The weakest condition known to ensure this is

\begin{defn}
\label{def:tempered-sequence}A F\o{}lner sequence $\left\{ F_{n}\right\} $
in $G$ is \emph{tempered} if for some constant $C$ independent of
$n$, $|(\bigcup _{k<n}F_{k}^{-1})F_{n}|\leq C|F_{n}|$ for every
$n$. 
\end{defn}
(We extend the group operations to sets in the usual manner, so $AB=\left\{ ab\, :\, a\in A\, ,\, b\in B\right\} $,
etc).

Throughout this paper, $\left\{ F_{n}\right\} $ will denote a tempered
F\o{}lner sequence which also satisfies $|F_{n}|\geq n$. It is not
hard to see that if $G$ is an infinite group any F\o{}lner sequence
in $G$ has a tempered subsequence satisfying this.

For a tempered F\o{}lner sequence, the ergodic theorem states that
for every $f\in L^{1}(\Omega ,\mathcal{F},P)$, \[
\lim _{n\rightarrow \infty }\frac{1}{|F_{n}|}\sum _{g\in F_{n}}f(g\omega )=\int f\, dP\]
 almost everywhere and in $L^{1}$. 

The entropy of a process $\{X_{h}\}_{g\in G}$ is defined as in the
case of $\mathbb{Z}$-processes by \[
h=\lim _{k\rightarrow \infty }\frac{1}{|F_{k}|}H(\{X_{g}\}_{g\in F_{k}})\]
where $H(\cdot )$ is the usual entropy function defined for any finite-valued
random variable. For a discussion of entropy theory of amenable processes
see \cite{W03}. 

For such processes, the Shannon-McMillan-Breiman (SMB) theorem states
that for $\{F_{n}\}$ a tempered sequence, \[
\frac{1}{|F_{n}|}\log P([x(F_{n})])\rightarrow h\]
 for almost every $x\in \Omega $. In this generality, these theorems
are due to E. Lindenstrauss \cite{Lin01}; see also \cite{W03} and
\cite{OW83}.

\subsection{Return times}

In the return-time theorems stated above for $\mathbb{Z}$ and $\mathbb{Z}^{d}$
processes a single F\o{}lner sequence $\{F_{n}\}$ (segments or cubes)
was used both to define the {}``shape'' of the patterns whose repetitions
we are looking for and as the shape of the domain in which we seek
these repetitions. There is no reason for the same sequence to play
both these roles, and it will be useful to separate them. In our formulation
of the problem we will consider two sequences of finite subsets of
$G$: $\{F_{n}\}$ ($F$ for F\o{}lner) which will serve as the shape
of the patterns, and $\{W_{n}\}$ ($W$ for window) which will serve
as the set within which we search for repetitions. Throughout, $\{F_{n}\}$
will be a tempered F\o{}lner sequence satisfying $|F_{n}|\geq n$.
The properties we demand of $\left\{ W_{n}\right\} $ will vary.

We will be interested in the recurrence of $F_{k}$-patterns in the
windows $W_{n}$. We say that $x(F_{k})$ \emph{repeats} at a point
$f\in G$ if $fx\in [x(F_{k})]$, ie if $fx(g)=x(g)$ for every $g\in F_{k}$.
If $x(F_{k})$ repeats at $f$, we say that the repetition is disjoint
from $F_{k}$ (or from $x(F_{k})$) if $F_{k}\cap F_{k}f=\emptyset $.

Define\[
R_{k}(x)=R_{k}^{(F,W)}(x)=\min \left\{ n\, \left|\begin{array}{c}
 \textrm{there exists a disjoint}\\
 \textrm{repetition of }x(F_{k})\textrm{ in }W_{n}\end{array}
\right.\right\} \]
 or $R_{k}(x)=\infty $ if the set on the right in empty.

\subsection*{Remarks}

\begin{enumerate}
\item $R_{k}(x)$ is the first {}``return index'' of $x$ to $[x(F_{k})]$,
\emph{excluding} recurrences which intersect the original pattern
$x(F_{k})$. In this respect our definition differs from the definition
for the case $G=\mathbb{Z}^{d}$ from \cite{OW02} described in the
previous section, where it is not required that the repetition be
disjoint from the original occurrence of the pattern. 
\item We require only that the {}``center'' of the repetition $f$ to
be in $W_{n}$. We do not require the whole of the repetition to be
in $W_{n}$, ie we do not require $F_{k}f\subseteq W_{n}$. Note that
in the known cases of cubes, $F_{n},W_{n}=[-n;n]^{d}\subseteq \mathbb{Z}^{d}$,
requiring that the repetition of $x(F_{k})$ be completely contained
in $W_{R_{k}}$ does not change the asymptotic behavior of $R_{k}$.
In general, however, the definitions are not equivalent. Indeed it
is not necessary that any translate of $F_{k}$ be contained in any
$W_{n}$.
\item In what follows we will be working with several different sequences
and then the notation $R_{k}^{(F,W)}$ will be used to emphasize the
dependence on $\{F_{n}\}$ and $\{W_{n}\}$. 
\item If $\{F_{n}\}$ is any sequence (not necessarily F\o{}lner) and $\{W_{n}\}$
is a tempered F\o{}lner sequence then the ergodic theorem implies
that $R_{k}^{(F,W)}(x)$ is finite almost surely. In other cases it
is not clear that $R_{k}<\infty $, and for general window sequence
this may not be the case. 
\end{enumerate}
$R_{k}$ is not a very {}``stable'' quantity, and even when finite
we cannot hope for its growth rate to provide any significant information
about the process, since it strongly reflects the indexing of the
sequence $\left\{ W_{n}\right\} $. For any $G$-process and choice
of $\{F_{n}\}$ and $\{W_{n}\}$, the value of $R_{k}$ may be arbitrarily
increased simply by allowing repetitions in the $W_{n}$'s, or decreased
by thinning the window sequence out. A better measure of return time
is the volume (=size) of the window set in which the repetition is
observed. We therefore define\[
T_{k}=T_{k}^{(F,W)}=|W_{R_{k}^{(F,W)}}|\]
 or $T_{k}=\infty $ if $R_{k}=\infty $, and set \[
T^{*}=T_{(F,W)}^{*}=\limsup _{k\rightarrow \infty }\frac{1}{|F_{k}|}\log T_{k}^{(F,W)}\]
 \[
T_{*}=T_{*}^{(F,W)}=\liminf _{k\rightarrow \infty }\frac{1}{|F_{k}|}\log T_{k}^{(F,W)}.\]
 The aforementioned result for $\mathbb{Z}^{d}$ can be given in terms
of $T_{k}$, since for $W_{n}=[-n;n]^{d}$, the size of $W_{R_{k}}$
is $(2R_{k}+1)^{d}$; the result then is equivalent to $T^{*}=T_{*}=h$. 

One half of the $\mathbb{Z}^{d}$ return-time theorem is true in the
following very general form, which provides new information even about
$\mathbb{Z}$ processes:

\begin{thm}
\label{thm:intro-lower-bound-for-T}Let $G$ be a countable amenable
group and $\left\{ X_{g}\right\} _{g\in G}$ a finite-valued ergodic
$G$-process. If $\{F_{n}\}$ is a tempered F\o{}lner sequence in
$G$ and $\{W_{n}\}$ an increasing sequence of finite subsets of
$G$, then for almost every $x$, \[
T_{*}^{(F,W)}(x)=\liminf _{k\rightarrow \infty }\frac{1}{|F_{k}|}\log T_{k}^{(F,W)}(x)\geq h.\]

\end{thm}
One might hope that the bound $T^{*}\leq h$ is also true; if it were,
we would have $\lim \frac{1}{|F_{k}|}\log T_{k}=h$. This need not
be true for general window sequences. The asymmetry between $T_{*}$
and $T^{*}$, which does not appear in the special cases examined
in \cite{OW02}, is not due to some shortcoming of our methods, but
is inherent in the problem. To begin with, in order for the upper
bound to hold, at the very least some growth condition must be imposed
on $\{W_{n}\}$. For example if we set $F_{n}=W_{n}=[0;2^{|F_{n-1}|}]\subseteq \mathbb{Z}$
then $|F_{n+1}|>2^{|F_{n}|}$ so $\frac{1}{|F_{n}|}\log T_{k}\geq 1$
is always true, even if $h<1$. We therefore cannot hope for a good
upper bound on the growth rate of $T_{k}$ if we do not restrict the
growth of the window sequence $\{W_{n}\}$.

Even if the sequence $\{W_{n}\}$ grows slowly, the inequality $T^{*}\leq h$
need not hold. For a simple example, consider the process of alternating
sequences of $0$s and $1$s with equal probabilities for observing
$0$ or $1$ at any given place, so the possible realizations are
$\ldots 0101010\ldots $ and $\ldots 1010101\ldots $. This process
is clearly ergodic and has entropy $0$. If we take $F_{n}=[-n;n]$
and $W_{n}$ to be sets containing only odd integers, eg $W_{n}=[-n,n]\cap (2\mathbb{Z}+1)$,
then $R_{k}=\infty $ for every realization. In this example $\left\{ W_{n}\right\} $
is of course not a F\o{}lner sequence; in section \ref{subsec:bad-window-sequences}
we show that the upper bound may fail even when $\left\{ W_{n}\right\} $
are {}``almost'' segments and possess extremely good properties
which are known to be sufficient for the ergodic and SMB theorems
to hold.

In order to state our results about $T^{*}$ we will require $\left\{ W_{n}\right\} $
to have some additional combinatorial properties, which will play
a role, in various forms, also in the density theorems stated below.

\begin{defn}
\label{def:incremental}Let $\{W_{n}\}$ be some sequence of finite
subsets of $G$. A sequence of right-translates $\{W_{n(i)}f_{i}\}_{i=1}^{I}$
($f_{i}\in G$) is called \emph{incremental} if 
\begin{enumerate}
\item $n(1)\geq n(2)\geq \ldots \geq n(I)$.
\item $f_{i}\notin \cup _{j<i}W_{n(j)}f_{j}$ for $1\leq i\leq I$.
\end{enumerate}
\end{defn}
{}

\begin{defn}
\label{def:filling}Let $\left\{ W_{n}\right\} $ be a sequence of
finite subsets of $G$ with $1_{G}\in W_{n}$. $\left\{ W_{n}\right\} $
is \emph{filling with constant $C$} (or $C$-filling) if for every
incremental sequence $\{W_{n(i)}f_{i}\}$ it holds that $|\cup W_{n(i)}f_{i}|\geq C\sum |W_{n(i)}f_{i}|$.
A set $W$ is individually $C$-filling if for every incremental sequence
$\{Wf_{i}\}$ it holds that $|\cup Wf_{i}|\geq C\sum |Wf_{i}|$. 
\end{defn}
The statement that $\left\{ W_{n}\right\} $ is filling means, informally,
this: suppose one tries to cover $G$ with translates of the $W_{n}$,
using a {}``greedy'' algorithm, that is, putting down some translate
of a $W_{n}$, then another (possibly smaller) $W_{k}$ over a point
not yet covered, etc. Note that we allow the different translates
to intersect. Then if you stop after finitely many steps, the size
of what has in fact been covered will be proportional to the total
of the sizes of the $W_{n}$'s used up until that stage.

A fixed set $W$ is filling with constant $\frac{1}{|W|}$ because
in an incremental sequence each point is covered at most $|W|$ times,
but generally sequences $\left\{ W_{n}\right\} $ containing infinitely
many different sets will not be filling, as demonstrated by the sequence
$\left\{ [0;n]\right\} _{n=1}^{\infty }$ in $\mathbb{Z}$ (though
not filling, this is a tempered F\o{}lner sequence). On the other
hand one can verify that $\left\{ [-n;n]\right\} _{n=1}^{\infty }$
is filling in $\mathbb{Z}$. Still another example is the sequence
$\left\{ [-n;n]\cap 2\mathbb{Z}\right\} $ in $\mathbb{Z}$ which
is filling but not a F\o{}lner sequence. 

From these examples it should be clear that the notion of being filling
is rather delicate, and not related directly to amenability. In particular,
note that while a F\o{}lner sequence is in some sense a large subset
of the group, a filling sequence can be small, eg contained in a proper
subgroup.

\begin{defn}
\label{def:quasi-filling}An increasing sequence $\left\{ W_{n}\right\} $
of finite subsets of $G$ with $1_{G}\in W_{n}$ is \emph{quasi-filling}
if each $W_{n}$ is $C_{n}$-filling and $\sum \frac{|W_{n}|^{-\alpha }}{C_{n}}<\infty $
for every $\alpha >0$. 
\end{defn}
Note that in general filling sequences need not be quasi filling.

\begin{thm}
\label{thm:intro-upper-bound-for-T}Let $G$ be a countable amenable
group and $\left\{ X_{g}\right\} _{g\in G}$ a finite-valued ergodic
$G$-process. If $\{F_{n}\}$ is a tempered F\o{}lner sequence in
$G$ and $\{W_{n}\}$ either a filling or a quasi-filling sequence
such that $\lim _{n}\frac{\log |W_{n}|}{\log |W_{n-1}|}=1$, then
\[
T^{*}=\limsup _{k\rightarrow \infty }\frac{1}{|F_{k}|}\log T_{k}^{(F,W)}\leq h\]
almost surely.
\end{thm}
\begin{cor}
If $\{F_{n}\},\{W_{n}\}$ satisfy the conditions of theorem \ref{thm:intro-upper-bound-for-T}
and \ref{thm:lower-bound-for-T}, then \[
\lim _{k\rightarrow \infty }\frac{1}{|F_{k}|}\log T_{k}^{(F,W)}=h\]
almost surely.
\end{cor}
It is interesting to note that in \cite{OW02}, the return times theorem
was proved also for the case $F_{n}=W_{n}=[0;n]^{d}\subseteq \mathbb{Z}^{d}$,
which is not covered by our results, because $\{[0;n]^{d}\}$ is not
(quasi) filling. We will not discuss this here.

\subsection{\label{subsec:recurrence-densities} Recurrence densities}

As we have seen, if $\{W_{n}\}$ grows too swiftly it may not be possible
for $T^{*}\leq h$ to hold, because the first $W_{n}$ within which
we observe $x(F_{k})$ may be too large compared to $F_{k}$. To compensate
for this, one expects that the number of times the pattern $x(F_{k})$
appears in $W_{n}$ will grow along with $|W_{n}|$.

Another way of looking at it is this: Suppose $\left\{ F_{n}\right\} ,\left\{ W_{n}\right\} $
are sequences such that the return times theorem of the last section
holds: $\frac{1}{|F_{k}|}\log T_{k}^{(F,W)}\rightarrow h$. If $\left\{ W_{n}\right\} $
grows slowly enough, then one sees from the return time theorem that
the density of recurrences of $x(F_{k})$ in $W_{R_{k}(x)-1}$ is
approximately $2^{-h|F_{k}|}$, since only the central copy of the
pattern exists. Now suppose in addition that $\left\{ W_{n}\right\} $
is itself a tempered F\o{}lner sequence; since by the SMB theorem
the probability of the atom $[x(F_{k})]$ is approximately $2^{-h|F_{k}|}$,
the ergodic theorem tells us that for large enough $n$ (depending
on $k$ and $x$) the frequency with which $x(F_{k})$ repeats in
$x(W_{n})$ is approximately equal to the probability of the atom
$[x(F_{k})]$. So we see that for large $k$, the density of recurrences
of $x(F_{k})$ in $x(W_{n})$ is close to $2^{-h|F_{k}|}$ both when
$n$ is small (but large enough for the question to be meaningful)
and for large enough $n$. However, neither the return-time theorem
nor the ergodic theorem give us any information about what happens
for intermediate values of $n$.

We define the recurrence frequency of $x(F_{k})$ in $x(F_{n})$ to
be \[
U_{k,n}(x)=U_{k,n}^{(F,W)}(x)=\frac{1}{|W_{n}|}\max \left\{ |E|\, \left|\begin{array}{c}
 1_{G}\in E\subseteq W_{n}\, ,\, \textrm{the collection}\\
 \{F_{k}f\}_{f\in E}\textrm{ is pairwise disjoint},\\
 \textrm{and }x(F_{k})\textrm{ repeats at }f\end{array}
\right.\right\} \]
 and\[
U^{*}=U^{(F,W)*}=\limsup _{k\rightarrow \infty }\, \sup _{n}-\frac{1}{|F_{k}|}\log U_{k,n}^{(F,W)}\]
 \[
U_{*}=U_{*}^{(F,W)}=\liminf _{k\rightarrow \infty }\, \inf _{n}-\frac{1}{|F_{k}|}\log (U_{k,n}^{(F,W)}-\frac{1}{|W_{n}|}).\]

\subsection*{Remarks:}

\begin{enumerate}
\item The quantity $U_{k,n}$ was defined as the (maximal) density of \emph{disjoint}
repetitions of $x(F_{k})$ in $W_{n}$, one of which is the original
pattern. One may attempt to define $U_{k,n}$ differently, eg by counting
every repetition or counting every repetition disjoint from the original
pattern (but not mutually disjoint). These variations are discussed
in more detail in sections \ref{subsec:upper-bound-for-U} and \ref{subsec:lower-bound-for-U}. 
\item The reason for the correction factor $-\frac{1}{|W_{n}|}$ in the
definition of $U_{*}$ is that our goal is to prove $U_{*}\geq h$,
but in our definition of $U_{k,n}$ we allowed the {}``central''
copy of $x(F_{k})$ to be counted. This being the case, we cannot
expect the observed frequency $U_{k,n}$ to drop below $2^{-h|F_{k}|}$
if the size $|W_{n}|$ is not on the order of $2^{h|F_{k}|}$. 
\end{enumerate}
As in the case of $T^{*}$, in order to prove $U_{*}=U^{*}=h$ we
need some assumptions about the combinatorial properties of the window
sequence $\{W_{n}\}$. For the bound $U^{*}\leq h$, it is enough
that $\{W_{n}\}$ be filling (definition \ref{def:filling} above):

\begin{thm}
\label{thm:intro-upper-bound-for-U}Let $G$ be a countable amenable
group and $\left\{ X_{g}\right\} _{g\in G}$ a finite-valued ergodic
$G$-process. Let $\{F_{n}\}$ be a tempered F\o{}lner sequence in
$G$. If $\{W_{n}\}$ is a filling sequence, then \[
U^{*}=\limsup _{k\rightarrow \infty }\, \sup _{n}-\frac{1}{|F_{k}|}\log U_{k,n}^{(F,W)}\leq h\]
almost surely.
\end{thm}
See also theorem \ref{thm:upper-bound-for-U-weak-version} below for
a slightly weaker version which applies to quasi-filling sequences.

For the lower bound $U_{*}\geq h$ we need an even stronger notion
of {}``filling up space'':

\begin{defn}
\label{def:incompressible}An increasing sequence $\left\{ W_{n}\right\} $
of finite subsets of $G$ with $1_{G}\in W_{n}$ is said to be \emph{incompressible
with constant $C$} (or $C$-incompressible) if for any incremental
sequence $\left\{ W_{n(i)}f_{i}\right\} $, the number of the sets
$W_{n(i)}f_{i}$ containing $1_{G}$ is at most $C$. A finite set
$W\subseteq G$ containing $1_{G}$ is \emph{individually $C$-incompressible}
if for any incremental sequence $\left\{ Wf_{i}\right\} $, the number
of the sets $Wf_{i}$ containing $1_{G}$ is at most $C$. 
\end{defn}
Clearly the condition that $1_{G}$ be in at most $C$ sets is equivalent
to requiring that every $g\in G$ is contained in at most $C$ sets.
It follows that a $C$-incompressible sequence is $\frac{1}{C}$-filling.

Any finite set $W$ containing $1_{G}$ is $|W|$-incompressible.
Ascending sequences of subgroups, and ascending centered cubes in
$\mathbb{Z}^{d}$, are examples of incompressible sequences (see also
section \ref{sec:examples}).

We also need the following:

\begin{defn}
\label{def:interpolation} Let $\{F_{n}\},\{W_{n}\}$ be sequences
of finite subsets of $G$. A sequence $\{Y_{n}\}$ of finite subsets
of $G$ is called an \emph{interpolation sequence} \emph{for} $\{F_{n}\},\{W_{n}\}$
if 
\begin{enumerate}
\item $\left\{ Y_{n}\right\} $ is filling.
\item There exists a constant $C$ such that if $|W_{n}|\geq |Y_{k}|$ then
$|Y_{k}W_{n}|\leq C|W_{n}|$.
\item For every pair of real numbers $0\leq \alpha <\beta $, for every
large enough $n$ there is an index $k$ such that $2^{\alpha |F_{n}|}\leq |Y_{k}|\leq 2^{\beta |F_{n}|}$.
\end{enumerate}
In the case where $F_{n}=W_{n}$, we say simply that $\left\{ Y_{k}\right\} $
is an interpolation sequence for $\left\{ F_{n}\right\} $. 
\end{defn}
Note that $\{F_{n}\}$ appears only in condition (3) of the definition,
and will be satisfied automatically (for any sequence $\left\{ F_{n}\right\} $)
if $\lim _{n}\frac{\log |Y_{n}|}{\log |Y_{n-1}|}=1$.

As an example, one may consider the sequence $\left\{ [-n;n]\right\} $
of intervals in $\mathbb{Z}$, which is its own interpolation sequence.
In fact if $G$ contains an element of infinite order, and we identify
the infinite-cyclic group it generates with $\mathbb{Z}$, then segments
in this subgroup can often serve as an interpolation sequence for
a suitably chosen window sequence in $G$. We come back to this in
section \ref{subsec:recurrence-in-subgroups}.

\begin{thm}
\label{thm:intro-lower-bound-for-U}Let $G$ be a countable amenable
group and $\left\{ X_{g}\right\} _{g\in G}$ a finite-valued ergodic
$G$-process. Let $\{F_{n}\}$ be a tempered F\o{}lner sequence in
$G$. Let $\{W_{n}\}$ be an incompressible sequence. If there exists
an interpolation sequence for $\{F_{n}\},\{W_{n}\}$, then \[
U_{*}^{(F,W)}=\liminf _{k\rightarrow \infty }\, \inf _{n}-\frac{1}{|F_{k}|}\log (U_{k,n}^{(F,W)}-\frac{1}{|W_{n}|})\geq h\]
almost surely.
\end{thm}
Since $T_{*}\geq h$ almost surely, we have $U_{*}=\liminf _{k}\, \inf _{n}-\frac{1}{|F_{k}|}\log U_{k,R_{k}}$
(ie the correction term may be dropped), so 

\begin{cor}
If $\left\{ F_{n}\right\} $,$\left\{ W_{n}\right\} $ are as in theorem
\ref{thm:intro-lower-bound-for-U}, then \[
\lim _{k}-\frac{1}{|F_{k}|}\log U_{k,R_{k}}^{(F,W)}=h\]
almost surely.
\end{cor}
A slightly weaker condition than incompressibility is the following:

\begin{defn}
\label{def:quasi-incompressible}Let $\left\{ W_{n}\right\} $ be
an increasing sequence of finite subsets of $G$ containing $1_{G}$,
and suppose that each $W_{n}$ is individually $C$-incompressible
(for the same $C$). $\{W_{n}\}$ is \emph{quasi-incompressible} if
for every $\lambda >0$ and large enough $k$, for every $n>k$ the
$W_{k}$ boundary of $W_{n}$ is at most a $|W_{k}|^{-\lambda }$
fraction of $W_{n}$, ie\[
\frac{1}{|W_{n}|}\{f\in W_{n}\, :\, W_{k}f\not \subseteq W_{n}\}\leq |W_{k}|^{-\lambda }.\]
 
\end{defn}
Note that with this definition, an incompressible sequence may not
be quasi-incompressible.

\begin{thm}
\label{thm:intro-alternate-lower-bound-for-U}Let $G$ be a countable
amenable group and $\left\{ X_{g}\right\} _{g\in G}$ a finite-valued
ergodic $G$-process. Let $\{F_{n}\}$ be a tempered F\o{}lner sequence
and $\{W_{n}\}$ a quasi-incompressible sequence in $G$. If there
exists an interpolation sequence for $\{F_{n}\},\{W_{n}\}$, then
for almost every $x$, \[
U_{*}^{(F,F)}(x)=\liminf _{k\rightarrow \infty }\, \inf _{n}-\frac{1}{|F_{k}|}\log (U_{k,n}^{(F,F)}(x)-\frac{1}{|F_{n}|})\geq h.\]

\end{thm}
\begin{rem*}
We have not required in any of the theorems that $\left\{ W_{n}\right\} $
be a F\o{}lner sequence. In fact one may verify that if $H<G$ and
$\left\{ W_{n}\right\} $ is a sequence of subsets of $H$ satisfying
one of the various filling properties described above (except for
the property of the existence of an interpolation sequence), then
the property continues to hold when $\left\{ W_{n}\right\} $ is viewed
as a sequence of subsets of $G$. This means that the recurrence phenomena
we have described can often be observed when one counts repetitions
in a subgroup of $G$. We will come back to this in more detail in
section \ref{subsec:recurrence-in-subgroups} below.
\end{rem*}
The rest of this paper is organized as follows. In section \ref{sec:examples}
we discuss various ways to construct good window sequences and classes
of groups in which there exist good window sequences, with respect
to which the recurrence phenomena hold. We also give an example showing
that other plausible combinatorial properties of $\left\{ W_{n}\right\} $
do not ensure that the upper bounds for $T^{*}$ and $U^{*}$ hold.
In section \ref{sec:upper-bounds} we prove the upper bounds $T^{*},U^{*}\leq h$.
In section \ref{sec:lower-bounds} we prove the lower bound $T_{*}\geq h$,
whose proof has a different flavor than the others, and $U_{*}\geq h$,
which depends strongly on the lower bound $T_{*}\geq h$. Finally
in section \ref{sec:open-questions} we state some open questions.

\begin{acknowledgement*}
The results in this paper are part of the author's MA thesis, conducted
under the guidance of Prof. B. Weiss, whom I would like to thank for
all his patience and good advice. 
\end{acknowledgement*}

\section{\label{sec:examples} Constructions, examples and counterexamples}

In this section we will examine the problem of the existence of good
window sets. We consider two questions: When do good window sets exist
at all, and when do {}``large'' window sets exist. By large we mean
window sequences which are themselves F\o{}lner sequences. We also
give an example which shows that, while the filling properties may
not be strictly necessary for our theorems to hold, other combinatorial
conditions on the window sequence, which are known to be sufficient
for other pointwise theorems, are not sufficient for upper bounds
on return times and densities.

\subsection{\label{subsec:general-constructions}General constructions}

Before going into concrete examples we examine how good window and
interpolation sequences can be constructed from other sequences.

\begin{lem}
\label{lemma:sequences-in-product-groups}Suppose $\{W'_{n}\}$ is
$C'$-incompressible in $G'$ and $\{W''_{n}\}$ is $C''$-incompressible
in $G''$. Set $W_{n}=W'_{n}\times W''_{n}$. The $\{W_{n}\}$ is
incompressible in $G=G'\times G''$. 
\end{lem}
\begin{proof}
Clearly $\{W_{n}\}$ is increasing and $1_{G}\in W_{n}$ for all $n$.
Suppose $\{W_{n(i)}g_{i}{}_{i=1}^{I}\}$ is incremental and $1_{G}\in \cap _{i}W_{n(i)}g_{i}$;
we must show that $I$ is bounded. Write $g_{i}=(g'_{i},g''_{i})$.
Since $1_{G'}\in \cap _{i}W'_{n(i)}g'_{i}$ and $1_{G''}\in \cap _{i}W''_{n(i)}g''_{i}$,
it suffices to show that there is a $C$ such that if $I>C$ then
either there is a subsequence $\left\{ i(1),\ldots ,i(C'+1)\right\} \subseteq \left\{ 1,\ldots ,I\right\} $
with $\{W'_{n(i(j))}g'_{i(j)}\}_{i=1}^{C'+1}$ incremental or there
is a subsequence $\left\{ i(1),\ldots ,i(C''+1)\right\} \subseteq \left\{ 1,\ldots ,I\right\} $
with $\{W'_{n(i(j))}g'_{i(j)}\}_{i=1}^{C''+1}$ incremental.

We prove this claim by induction on $C'+C''$. To be precise, we claim
that for each pair $C',C''$ there is a number $\alpha (C',C'')$
such that if $\{W_{n(i)}g_{i}\}_{i=1}^{I}$ is an incremental sequence
in $G$ of length $I>\alpha (C',C'')$ then either its projection
on the first coordinate has an incremental subsequence of length greater
than $C'$ or its projection on the second coordinate has an incremental
subsequence of length $C''$. Clearly \[
\alpha (C',1)=C'\; \alpha (1,C'')=C''\; \alpha (1,1)=1.\]
For the induction step, let $\left\{ W_{n(i)}g_{i}\right\} _{i=1}^{I}$
be incremental in $G$. Write \[
U_{1}=\left\{ 1<i\leq I\, :\, g'_{i}\notin W'_{n(1)}g'_{1}\right\} \; ,\; U_{2}=\left\{ 1<i\leq I\, :\, g''_{i}\notin W''_{n(1)}g''_{1}\right\} \]
 clearly $U_{1}\cup U_{2}=\left\{ 2,\ldots ,I\right\} $. If $|U_{1}|>\alpha (C'-1,C'')$
or $|U_{2}|>\alpha (C',C''-1)$ then $\left\{ W_{n(i)}g_{i}\right\} $
has a subsequence of the type we desire. But one of these relations
will certainly hold if $|U_{1}\cup U_{2}|>\alpha (C'-1,C'')+\alpha (C',C''-1)$,
so we can set \[
\alpha (C',C'')=\alpha (C'-1,C'')+\alpha (C',C''-1)+1.\qedhere \]

\end{proof}
The point here is that we obtain a {}``large'' incompressible sequence
in $G'\times G''$. {}``Small'' sequences exist since an incompressible
sequence in $G'$ (or $G''$) is incompressible in $G'\times G''$
with the natural embedding. Note that if $\left\{ W'_{n}\right\} ,\left\{ W''_{n}\right\} $
are F\o{}lner sequences in $G',G''$ respectively then $\left\{ W'_{n}\times W''_{n}\right\} $
is F\o{}lner in $G'\times G''$.

In \cite{W01} it was shown how certain properties of F\o{}lner sequences
can be pulled up through exact sequences. We rematk that we do not
know whether the existence of filling or incompressible F\o{}lner
sequences in $K$ and $G/K$ for a normal subgroup $K$ of $G$ implies
that such sequences exist in $G$, although it seems unlikely that
this is the case (again, the sequence in $K$ has the same good filling
properties in $G$, but is not F\o{}lner). 

The requirement that $\left\{ W_{n}\right\} $ be $C$-filling (or
$C$-incompressible) is a property of the entire sequence, and in
particular implies that each $W_{n}$ is $C$-filling ($C$-incompressible)
as an individual set. The existence of a constant $C$ such that each
$W_{n}$ is $C$-filling ($C$-incompressible) does not imply that
the same about the sequence. One can, however, sometimes pass to a
subsequence of the $\left\{ W_{n}\right\} $ which will be filling. 

\begin{lem}
\label{lemma:creating-filling-sequences}Suppose $\left\{ W_{n}\right\} $
is an increasing sequence of finite subsets of $G$ with $1_{G}\in W_{n}$
such that 
\begin{enumerate}
\item $\left\{ W_{n}\right\} $ is a F\o{}lner sequence.
\item For some constant $C$, each $W_{n}$ is individually $C$-filling.
\end{enumerate}
Then there exists a subsequence $\left\{ W_{n(i)}\right\} $ which
is filling. 
\end{lem}
\begin{proof}
Since $\{W_{n}\}$ is a F\o{}lner sequence, we may select a subsequence
$\{W_{n(i)}\}$ with the property that \[
\frac{1}{|W_{n(i)}|}\left|\left\{ g\in W_{n(i)}\, :\, W_{n(i-1)}g\subseteq W_{n(i)}\right\} \right|\geq \frac{C}{2}.\]
 Write $W'_{i}=W_{n(i)}$; we claim that $\{W'_{i}\}$ is filling.
Suppose that $\{W'_{i(j)}f_{j}\}_{j=1}^{J}$ is an incremental sequence.
Let $I(1)>\ldots >I(N)$ be all the values of $i(1),\ldots ,i(J)$,
and let $J_{k}=\left\{ 1\leq j\leq J\, :\, I(j)=k\right\} $. Since
$W_{I(1)}$ is $C$-filling, we have \[
|\bigcup _{j\in J_{1}}W_{I(1)}f_{j}|\geq C\sum _{j\in J_{1}}|W_{I_{1}}f_{j}|.\]
 Let \[
E_{1}=\{g\in \bigcup _{j\in J_{1}}W_{I(1)}f_{j}\, :\, W_{I(2)}g\subseteq \bigcup _{j\in J_{1}}W_{I(1)}f_{j}\}\]
 then \begin{eqnarray*}
|E_{1}| & \geq  & |\bigcup _{j\in J_{1}}W_{I(1)}f_{j}|-|\{g\in \bigcup _{j\in J_{1}}W_{I(1)}f_{j}\, :\, W_{I(2)}g\not \subseteq \bigcup _{j\in J_{1}}W_{I(1)}f_{j}\}|\\
 & \geq  & C\sum _{j\in J_{1}}|W_{I(1)}f_{j}|-\sum _{j\in J_{1}}\frac{C}{2}|W_{I(1)}f_{j}|\\
 & = & \frac{C}{2}\sum _{j\in J_{1}}|W_{I(1)}f_{j}|
\end{eqnarray*}
 so \[
|\bigcup _{j\in J_{1}}W_{I(1)}f_{j}\setminus \bigcup _{j\in J_{2}\cup \ldots \cup J_{N}}W_{i(j)}f_{j}|\geq |E_{1}|\geq \frac{C}{2}\sum _{j\in J_{1}}|W_{I_{1}}f_{j}|\]
 which means that the fraction of $\bigcup _{j\in J_{1}}W_{I(1)}f_{j}$
which is not contained in any of the {}``lower levels'' is at least
$\frac{C}{2}\sum _{j\in J_{1}}|W_{I(1)}f_{j}|$. Repeating this argument
for $I(1),\ldots ,I(N)$ we obtain\[
|\bigcup _{1\leq j\leq J}W_{i(j)}f_{j}|\geq \frac{C}{2}\sum _{1\leq j\leq J}|W_{i(j)}f_{j}|.\qedhere \]

\end{proof}
The following is proved similarly:

\begin{lem}
\label{lemma:creating-quasi-incompressible-sequences}Suppose $\left\{ W_{n}\right\} $
is an increasing sequence of finite subsets of $G$ with $1_{G}\in W_{n}$
such that 
\begin{enumerate}
\item $\left\{ W_{n}\right\} $ is a F\o{}lner sequence.
\item For some constant $C$, each $W_{n}$ is $C$-incompressible. 
\end{enumerate}
Then there exists a subsequence $\{W_{n(i)}\}$ which is quasi-incompressible
.
\end{lem}
The last construction in this section deals with the existence of
interpolation sequences.

\begin{lem}
\label{lemma:creating-interpolation-sequences}Let $\{F_{n}\}$ be
an increasing tempered F\o{}lner sequence. Suppose $\{Y_{n}\}$ is
increasing and filling. Then there exists a subsequence $\{F_{n(i)}\}$
of $\{F_{n}\}$ such that the sequence $Y_{n}$ is an interpolation
sequence for $\{F_{n(i)}\}$. 
\end{lem}
For the sequence $\{Y_{n}\}$ to be useful, $\{F_{n}\}$ must fulfill
the hypothesis of theorem \ref{thm:intro-lower-bound-for-U} and so
$\{F_{n}\}$ must also be incompressible or quasi-incompressible;
but this is not necessary for the proof of the lemma.

\begin{proof}
Let $\{[\alpha _{i},\beta _{i}]\}_{i=1}^{\infty }$ be an enumeration
of the rational intervals, $0\leq \alpha _{i}<\beta _{i}$. We select
$F_{n(k)}$ and sets $Z_{i,k}=Y_{m(i,k)}$ ($1\leq i<k$) inductively
in $k$, as follows: Given $k$ and assuming that we have defined
$F_{n(j)}$ and $Z_{i,j}$ for every $1\leq i<j<k$, we select $F_{n(k)}$
and $Z_{i,k}$ for $1\leq i<k$ so that

\renewcommand{\theenumi}{\alph{enumi}}
\begin{enumerate}
\item $Z_{i,k}$ satisfies $2^{\alpha _{i}|F_{n(k-1)}|}\leq |Z_{i,k}|\leq 2^{\beta _{i}|F_{n(k-1)}|}$.
\item $F_{n(k)}$ is large enough that for all $i<k$ it holds that $|F_{n(k)}|<2^{\alpha _{i}|F_{n(k)}|}$.
\item $|(\bigcup _{i<k-1}Z_{i,k-1})F_{n(k)}|\leq 2|F_{n(k)}|$.
\item $n(k)$ is large enough so that for $1\leq i\leq k$ and every $n>n(k)$
there exists $m$ such that $2^{\alpha _{i}|F_{n}|}<|Y_{m}|<2^{\beta _{i}|F_{n}|}$. 
\end{enumerate}
\renewcommand{\theenumi}{\arabic{enumi}}

One then verifies that $\{Z_{i,k}\}$ is an interpolation sequence
for $\{F_{n(k)}\}$. Briefly, (a) ensures that (3) of the definition
of interpolation sequences holds, while (b) and (c) ensure condition
(2), and (d) makes it possible to continue the construction ((a) is
possible for $k$ because for $j<k$ we chose $W_{n(j)}$ to satisfy
(d) ). 
\end{proof}

\subsection{\label{subsec:Z-and-Z-d}The groups $\mathbb{Z}$ and $\mathbb{Z}^{d}$}

In the case of $\mathbb{Z}$, any increasing sequence of symmetric
intervals $\left\{ I_{n}\right\} $ is $2$-incompressible. This is
simple to verify.

Next, $\mathbb{Z}^{d}$ is the product of $\mathbb{Z}$ with itself
$d$ times, we know from lemma \ref{lemma:sequences-in-product-groups}
that incompressible sequences exist. In fact,

\begin{prop}
If $\{I_{n}^{(i)}\}_{n=1}^{\infty }$ are increasing sequences of
symmetric intervals, $i=1,\ldots ,d$, then the sets $W_{n}=I_{n}^{(1)}\times \ldots \times I_{n}^{(d)}$
form an $2^{d}$-incompressible sequence. 
\end{prop}
Lemma \ref{lemma:sequences-in-product-groups} gives a poorer estimate
for the constant. That the constant $2^{d}$ is correct (and optimal)
for such sequences can be established directly.

We will call sets of the form in the proposition boxes\emph{.} If
$\left\{ W_{n}\right\} $ is an increasing symmetric sequence of boxes,
then it is easy to see that for any $n$ we can find boxes $W_{n}=Y_{n,0}\subseteq \ldots \subseteq Y_{n,k(n)}=W_{n+1}$
such that $|Y_{n,i+1}|/|Y_{n,i}|\leq 2$. Ordering all the $Y_{n,i}$
in a single increasing sequence $\left\{ Y_{m}\right\} $, we obtain
a filling sequence which grows at most exponentially. It is easy to
verify that $|W_{n}+Y_{m}|\leq 2^{d}|W_{n}|$ whenever $|W_{n}|\geq |Y_{m}|$;
so $\left\{ Y_{m}\right\} $ is an interpolation sequence for $\left\{ W_{n}\right\} $.
This gives

\begin{prop}
\label{prop:Z}For $G=\mathbb{Z}^{d}$, any increasing sequence $\left\{ W_{n}\right\} $
of boxes and any tempered F\o{}lner sequence $\left\{ F_{n}\right\} $
(and in particular $F_{n}=W_{n}$), $T_{*}^{(F,W)}\geq h$ and $U_{*}^{(F,W)}=U_{(F,W)}^{*}=h$.
If $\lim _{n}\frac{\log |W_{n+1}|}{\log |W_{n}|}=1$ then $T_{*}^{(F,W)}=T_{(F,W)}^{*}=h$. 
\end{prop}
This is essentially a reformulation of the results from \cite{OW02}.

One can show that certain increasing sequences of symmetric convex
sets in $\mathbb{Z}^{d}$ are incompressible as well.

\subsection{\label{subsec:Z-infinity}The group \protect$\mathbb{Z}^{\infty }$}

Consider the group $\mathbb{Z}^{\infty }=\oplus _{n=1}^{\infty }\mathbb{Z}$.
Identify the subgroups $\mathbb{Z}^{d}\times \left\{ (0,0\ldots )\right\} \subseteq \mathbb{Z}^{\infty }$
with $\mathbb{Z}^{d}$ in the obvious way.

\begin{prop}
Let $A\subseteq \mathbb{Z}^{d}$ have full dimension (ie $\left\langle A\right\rangle $
is of abelian rank $d$). If $A$ is $C$-incompressible then $C>d$. 
\end{prop}
\begin{proof}
For a set $A\subseteq \mathbb{Q}^{d}$, we say that $u\in \mathbb{Q}^{d}$
is an extreme point of $A$ if it is an extreme point of $\textrm{conv}_{\mathbb{Q}}A$.

We prove the following: If $A\subseteq \mathbb{Q}^{d}$ is finite
and has $\dim \left\langle A\right\rangle =d$, and $0$ is an extreme
point of $A$, then there exists a set $\left\{ u_{1},\ldots ,u_{d+1}\right\} \subseteq -A$
such that $\left\{ A+u_{i}\right\} _{i=1}^{d+1}$ is incremental.
In this case, clearly $0\in \cap (A+u_{i})$, so $A$ cannot be $C$-incompressible
for $C\leq d$.

The proof is by induction on $d$. For $d=1$ the claim is trivial:
take $a_{1}=0$ and $a_{2}\in A\setminus \{0\}$.

Now suppose we have proved it for $d-1$. Let $A\subseteq \mathbb{Q}^{d}$
be as in the claim. Set $E=\textrm{ext}\, \textrm{conv}\, A$, so
by assumption $0\in E$. Since $A$ is $d$-dimensional, so is $E$.
Select a linearly independent set $a_{1},\ldots ,a_{d-1}\in E\setminus \left\{ 0\right\} $.
Write $V=\textrm{span}_{\mathbb{Q}}\left\{ a_{1},\ldots ,a_{d-1}\right\} $,
and assume $\left\{ a_{i}\right\} $ was selected so that $A$ is
on one side of the subspace $V$, ie there is a linear functional
$\Lambda $ with $\Lambda (A)\subseteq \mathbb{Q}^{+}$ and $V=\ker \Lambda $.

Consider $A'=A\cap V$. This is a finite set of full dimension in
$V$ containing $0$ as an extreme point. Therefore there exist $u_{1},\ldots ,u_{d}\in -A'$
with $\left\{ A'+u_{i}\right\} _{i=1}^{d}$ incremental. Clearly $\left\{ A+u_{i}\right\} _{i=1}^{d}$
is also incremental. Now we need only note that for any $u\in A\setminus A'$,
we have that $-u\notin \cup _{i=1}^{d}A+u_{i}$ because $\Lambda (u)<0$
whereas $\Lambda (a+u_{i})=\Lambda (a)+0\geq 0$ for all $a\in A$
and $i=1,\ldots ,d$. Thus $\left\{ u_{1},u_{2},\ldots ,u_{d},-u\right\} $
is the set we are looking for. This completes the induction step.

The proposition now follows from the fact that it holds for $A$ iff
it holds for every translate of $A$, and we can always translate
$A$ so that $0$ is an extreme point. 
\end{proof}
\begin{cor}
There are no incompressible F\o{}lner sequences in $\mathbb{Z}^{\infty }$.
\end{cor}
\begin{proof}
If $\left\{ F_{n}\right\} $ is a F\o{}lner sequence in $\mathbb{Z}^{\infty }$
then for any $d$, $F_{n}$ must be eventually $d$-dimensional. 
\end{proof}
It seems likely that there is a similar bound on how filling a $d$-dimensional
set in $\mathbb{Z}^{d}$ can be, and this would imply that there are
no filling F\o{}lner sequences in $\mathbb{Z}^{\infty }$, but we
do not have a proof.

On the bright side, there do exist quasi-filling F\o{}lner sequences
in $\mathbb{Z}^{\infty }$, as we now demonstrate.

Suppose $I_{n}^{(i)}\subseteq \mathbb{Z}$ are symmetric segments
and $W_{n}=I_{n}^{(1)}\times \ldots \times I_{n}^{(d(n))}\subseteq \mathbb{Z}^{\infty }$
is an increasing sequence of finite-dimensional boxes. We know that
each set $W_{n}$ is $2^{-d(n)}$-filling in $\mathbb{Z}^{d(n)}$
and it follows that the same is true as subsets of $\mathbb{Z}^{\infty }$.
In order for $\left\{ W_{n}\right\} $ to be quasi filling, we must
have \[
\sum _{n}2^{d(n)}\cdot |W_{n}|^{-\alpha }<\infty \]
 for every $\alpha >0$. For this it is enough that $\left\{ W_{n}\right\} $
grow exponentially and $d(n)=O(\log |W_{n}|)$. We would also like
$\left\{ W_{n}\right\} $ to be a tempered F\o{}lner sequence. Such
a sequence can easily be constructed. For example, we may take \[
W_{n}=[-2^{n^{2}};2^{n^{2}}]^{[\log n]}.\]
 $\left\{ W_{n}\right\} $ is then an increasing quasi-filling F\o{}lner
sequence. For temperedness it suffices to check that $|W_{n-1}^{-1}W_{n}|\leq C|W_{n}|$.
We have\[
W_{n-1}^{-1}W_{n}=[-2^{n^{2}}-2^{(n-1)^{2}},2^{n^{2}}+2^{(n-1)^{2}}]^{[\log n]}\]
 so\[
\frac{|W_{n-1}^{-1}W_{n}|}{|W_{n}|}\leq \left(\frac{2\cdot (2^{n^{2}}+\cdot 2^{(n-1)^{2}})+1}{2\cdot 2^{n^{2}}+1}\right)^{[\log n]}\leq \]
 \[
\leq \left(1+2^{-2n+2}\right)^{[\log n]}\rightarrow 1.\]
 One also readily verifies that $\left\{ W_{n}\right\} $ satisfies
$\frac{\log |W_{n+1}|}{\log |W_{n}|}\rightarrow 1$. We therefore
have

\begin{prop}
\label{prop:Z-infinity}There exists a tempered F\o{}lner sequence
$\left\{ W_{n}\right\} $ in $\mathbb{Z}^{\infty }$ such that for
any tempered sequence $\left\{ F_{n}\right\} $ (and in particular
$F_{n}=W_{n}$) we have $T_{*}^{(F,W)}=T_{(F,W)}^{*}=h$. 
\end{prop}

\subsection{\label{subsec:locally-finite-groups} Locally finite groups}

$G$ is \emph{locally finite} if every finitely generated subgroup
is finite. Since we only consider countable groups this is equivalent
to saying that there are finite subgroups $G_{1}<G_{2}<\ldots $ whose
union is all of $G$.

It is easy to see that such a sequence $\{G_{n}\}_{n=1}^{\infty }$
is $1$-incompressible (and hence $1$-filling). Translates of $G_{i}$'s
are just cosets, so if $f\notin G_{i}g$ then $G_{i}f\cap G_{i}g=\emptyset $.
Thus if $\{G_{n(i)}g_{i}\}_{i=1}^{I}$ is incremental, the $G_{n(i)}g_{i}$
are pairwise disjoint, so every $f\in G$ belongs to at most one translate.
Such a sequence $\{G_{n}\}$ is also clearly a tempered F\o{}lner
sequence, since every $g\in G$ is eventually a member of $G_{n}$
for large enough $n$.

If in addition $\{G_{i}\}$ grows slowly enough, then $\left\{ G_{i}\right\} $
is its own interpolation sequence. Thus

\begin{prop}
\label{prop:locally-finite-groups-have-good-sequences}If $G_{1}<G_{2}<\ldots $
and $G=\cup G_{i}$, then for $F_{k}=W_{k}=G_{k}$ we have that $T_{*}\geq h$
and $U^{*}\leq h$. If in addition $\frac{\log |G_{i+1}|}{\log |G_{i}|}\rightarrow 1$
then $T_{*}=T^{*}=h$ and $U_{*}=U^{*}=h$. 
\end{prop}
In the case where $\left\{ G_{i}\right\} $ grows too quickly to be
its own interpolation sequence, we can still show that for a suitably
chosen subsequence $\left\{ G_{i(k)}\right\} $ there exists an interpolation
sequence. This is based on

\begin{thm*}
(P. Hall and C. R. Kulatilaka, \cite{HK64}) Every infinite locally
finite group has an infinite abelian subgroup 
\end{thm*}
To use the theorem we need

\begin{prop}
An infinite locally finite abelian group $A$ possesses an increasing
sequence of finite subsets $(Y_{n})_{n=1}^{\infty }$each of which
is indivdually $4$-incompressible and $|Y_{n-1}|<|Y_{n}|\leq 2|Y_{n-1}|$.
\end{prop}
\begin{proof}
We can find finite subgroups $1<\ldots <A_{n}<\ldots A$ such that
$A_{n+1}/A_{n}$ is cyclic. Thus $A_{n+1}=A'_{n+1}\oplus C_{n+1}$
and $A_{n}=A'_{n+1}\oplus k_{n}C_{n+1}$ for $C_{n+1}$ some cyclic
group and $k_{n}\in \mathbb{N}$. Let $Z_{n,i}$ be a $2$-incompressible
sequence in $C_{n+1}$; then $Y_{n,i}=A_{n}\cup (A'_{n+1}\times Z_{n,i})$
is $4$-incompressible, and ordering all the $Y_{n,i}$ by inclusion
gives the desired sequence.
\end{proof}
Now if $G$ is locally finite, and $A<G$ an infinite abelian group,
we can find a slowly-growing sequence $(Y_{n})$ of finite subsets
of $A$ each of which is $4$-incompressible. Thus $(Y_{n})$ enjoys
the same properties as subsets of $G$, and we may apply lemma \ref{lemma:creating-interpolation-sequences}
to obtain

\begin{prop}
For any locally finite group $G=\cup _{n=1}^{\infty }G_{n}$ there
exists a sequence $G_{n(i)}$ such that for $F_{i}=G_{n(i)}$, \[
U_{*}^{(F,F)}=U_{(F,F)}^{*}=h.\]

\end{prop}
We remark that in order to obtain an interpolation sequence in $G$
one does not need the full force of the Hall-Kulatilaka theorem, but
can rather use the fact that for every $N$ there is a $K$ such that
every finite group $H$ with $|H|>K$ has an abelian subgroup of size
at least $N$, and therefore infinite locally finite groups have arbitrarily
large abelian subgroups. This can be proved by elementary methods.

\subsection{\label{subsec:groups-with-polynomial-growth} Groups with polynomial
growth}

Let $G$ be a finitely generated group and $\Gamma \subseteq G$ a
finite set of generators. For convenience we assume throughout that
$\Gamma $ is symmetric, ie that $\Gamma =\Gamma ^{-1}$. For every
element $g\in G$ define the \emph{length} of $g$ (with respect to
$\Gamma $) by\[
\ell _{\Gamma }(g)=\min \{k\, :\, g\in \Gamma ^{k}\}\]
 so $\ell _{\Gamma }(g)$ is the least number of elements of $\Gamma $
which need be multiplied together to give $g$.

It is easy to check that \[
d_{\Gamma }(g_{1},g_{2})=\ell _{\Gamma }(g_{1}^{-1}g_{2})\]
 defines a left invariant metric on $G$ (this is just the shortest
distance metric on $G$'s Cayley graph). The ball of radius $n$ in
this metric centered at the unit element of $G$ is \[
B_{n}^{\Gamma }=\{g\in G\, :\, \ell _{\Gamma }(g)\leq n\}.\]
 The \emph{growth function} of $G$ (with respect to $\Gamma $) is
defined by \[
\gamma _{_{\Gamma }}(n)=|B_{n}^{\Gamma }|=\#\{g\in G\, :\, \ell _{\Gamma }(g)\leq n\}.\]

A finitely generated group $G$ is said to have \emph{polynomial growth}
if its growth function $\gamma _{_{\Gamma }}$ is bounded by a polynomial,
and \emph{subexponential growth} if $\frac{1}{n}\log \gamma _{_{\Gamma }}(n)\rightarrow 0$.
The property of having polynomial/subexponential growth is really
a property of the group and does not depend on the particular generating
set $\Gamma $. To see this, let $\Gamma _{1},\Gamma _{2}$ be two
finite generating sets for $G$ as above. Then there is an $N$ for
which $\Gamma _{1}\subseteq (\Gamma _{2})^{N}$, so $(\Gamma _{1})^{n}\subseteq (\Gamma _{2})^{Nn}$.
This means that\[
\ell _{\Gamma _{1}}(g)\leq N\ell _{\Gamma _{2}}(g)\]
 and so\[
\gamma _{_{\Gamma _{1}}}(n)\leq \gamma _{_{\Gamma _{2}}}(Nn).\]
 From this, and the fact that the roles of $\Gamma _{1}$ and $\Gamma _{2}$
are interchangeable, it follows that $\gamma _{_{\Gamma _{1}}}$ is
bounded from above by a polynomial $p(x)$ of degree $d$ iff $\gamma _{_{\Gamma _{2}}}$
is also (though not necessarily the same polynomial), and similarly
$\frac{1}{n}\log \gamma _{_{\Gamma _{1}}}(n)\rightarrow 0$ iff $\frac{1}{n}\log \gamma _{_{\Gamma _{2}}}(n)\rightarrow 0$.

Though formally finite groups have polynomial growth, we will assume
from here on that our groups are infinite. We also assume that $\Gamma $
is a fixed finite symmetric generating set for $G$ and write $\gamma $
for $\gamma _{_{\Gamma }}$ and $B_{n}$ for $B_{n}^{\Gamma }$.

One can show using elementary methods that groups with subexponential
growth are amenable. Indeed,

\begin{prop}
There is a subset $I\subseteq \mathbb{N}$ of density $1$ such that
$\{B_{i}\}_{i\in I}$ is a F\o{}lner sequence.
\end{prop}
\begin{proof}
For suppose that $\gamma $ grows subexponentially. We will find a
subsequence $\{B_{n_{k}}\}_{k=1}^{\infty }$ of balls which is a F\o{}lner
sequence.

Write $\alpha (n)=\frac{\gamma (n+1)}{\gamma (n)}$ and for $\varepsilon >0$
set $I_{\varepsilon }=\{n\in \mathbb{N}\, :\, \alpha (n)\geq 1+\varepsilon \}$.
Since $\gamma $ grows subexponentially, for every $\varepsilon >0$,
$I_{\varepsilon }$ must have zero density; ie\[
\lim _{n\rightarrow \infty }\frac{1}{n+1}|I_{\varepsilon }\cap [0;n]|=0.\]
 Choose $N_{k}$ growing very rapidly such that \[
\frac{1}{N_{k}+1}|I_{1/k}\cap [0;N_{k}]|<1/k.\]
And define $I=\cup _{k}(I_{1/k}\cap [N_{k-1};N_{k}])$. Since if $n$
is such that $\alpha (n)\geq 1-\frac{1}{k}$ then $gB_{n}\subseteq B_{n+1}$
and \[
|gB_{n}\Delta B_{n}|\leq |B_{n+1}\setminus gB_{n}|+|B_{n+1}\setminus B_{n}|\leq \frac{2}{k}|B_{n}|.\]
We see that $(B_{n})_{n\in I}$ is F\o{}lner, and clearly $I$ has
density $1$.
\end{proof}
We will need the following stronger fact:

\begin{thm*}
(\cite{Ba72} and \cite{Grom81}) If $G=\left\langle \Gamma \right\rangle $
has polynomial growth then there is an integer $d$ and constants
$c_{1},c_{2}$such that $c_{1}n^{d}\leq \gamma _{\Gamma }(n)\leq c_{2}n^{d}$.
\end{thm*}
This implies that \[
|B_{n}^{-1}B_{n}|=|B_{n}^{2}|\leq c_{2}(2n)^{d}\leq \frac{2^{d}c_{2}}{c_{1}}|B_{n}|\]
so the F\o{}lner sequence from the proposition is in fact tempered.

We also conclude that every ball $B_{n}$ is $4^{d}c_{2}/c_{1}$ incompressible.
For suppose that $\{B_{n}g_{i}\}_{i=1}^{N}$ is incremental and $1_{G}\in \cap _{i}B_{n}g_{i}$.
Then $\{B_{n/2}g_{i}\}$ is a disjoint collection contained in $B_{2n}g_{i}$
so \[
N\cdot c_{1}(\frac{n}{2})^{d}\leq N\cdot \gamma _{n/2}=\sum _{i=1}^{N}\gamma _{n/2}\leq \gamma _{2n}\leq c_{2}(2n)^{d}\]
 and therefore $N\leq \frac{c_{2}}{c_{1}}4^{d}$. 

We thus get

\begin{prop}
If $G$ has polynomial growth there is a F\o{}lner sequence $\{F_{n}\}$
such that $\frac{1}{|F_{n}|}\log T_{k}^{(F,F)}\rightarrow h$ almost
surely.
\end{prop}
Using the construction from \ref{lemma:creating-filling-sequences}
and \ref{lemma:creating-interpolation-sequences} we also get

\begin{prop}
If $G$ has polynomial growth there is a F\o{}lner sequence $\{F_{n}\}$
such that $U_{(F,F)}^{*}=U_{*}^{(F,F)}=h$ almost surely.
\end{prop}
We also note that by Gromov's theorem \cite{Grom81} every infinite
group $G$ with polynomial growth has an element of inifinite order,
so the next section applies.

\subsection{\label{subsec:recurrence-in-subgroups} Recurrence in subgroups}

It is simple to verify that all the various filling-type properties
described so far are preserved in moving from subgroup to supergroup.
This immediately gives

\begin{prop}
\label{prop:filling-persists-from-subgroups-to-supergroup}If $H<G$
and $\{W_{n}\}$ is a filling or quasi-filling sequence in $H$ then
for any F\o{}lner sequence $\{F_{n}\}$ in $G$, $U_{(F,W)}^{*}\leq h$.
If in addition $\frac{\log |W_{n+1}|}{\log |W_{n}|}\rightarrow 1$
the $T_{*}^{(F,W)}=T_{(F,W)}^{*}=h$. 
\end{prop}
To obtain the same result for $U_{*}$, more care is needed, because
one cannot automatically transfer an interpolation sequence from a
subgroup to supergroup; this is due to the fact that condition (3)
of the definition of interpolation sequences (definition \ref{def:interpolation})
depends on the F\o{}lner sequence involved, and not only on the window
sequence, and the F\o{}lner sequences in the supergroup are quite
different from those in the subgroup. However, if $\{Y_{n}\}$ satisfies
(1) and (2) of the definition of interpolation sequences with respect
some to $\left\{ W_{n}\right\} $ and in addition $\lim _{n\rightarrow \infty }\frac{\log |Y_{n}|}{\log |Y_{n-1}|}=1$,
then (3) of the definition of interpolation sequences is satisfied
for \emph{any} sequence $\{F_{n}\}$, and so $\{Y_{n}\}$ is an interpolation
sequence in $G$ for $\{F_{n}\},\{W_{n}\}$ for \emph{any} sequence
$\{F_{n}\}$ in $G$. This gives

\begin{prop}
\label{prop:incompressibility-persists-from-subgroup-to-supergroup}If
$H<G$ and $\{W_{n}\}$ is an increasing incompressible sequence in
$H$, and there exists a sequence $\{Y_{n}\}$ satisfying (1) and
(2) of definition \ref{def:interpolation} and such that $\frac{\log |W_{n+1}|}{\log |W_{n}|}\rightarrow 1$,
then for any F\o{}lner sequence $\left\{ F_{n}\right\} $ in $G$
we have $U^{(F,W)*}=U_{*}^{(F,W)}=h$. 
\end{prop}
\begin{cor}
If $G$ contains an element of infinite order then there exist in
$G$ window sequences $\left\{ W_{n}\right\} $ for which for any
tempered F\o{}lner sequence $\left\{ F_{n}\right\} $, $U_{*}^{(F,W)}=U^{(F,W)*}=T_{*}^{(F,W)}=T^{(F,W)*}=h$.
\end{cor}
\begin{proof}
If $\left\langle g\right\rangle \cong \mathbb{Z}$ take $W_{n}=[-n;n]\subseteq \left\langle g\right\rangle $;
since $\left\{ W_{n}\right\} $ is incompressible in $\mathbb{Z}$
the same is true in $G$, and it is its own interpolation sequence.
\end{proof}
We remark that in the case that $\{W_{n}\}$ is contained in some
proper subgroup $H$ of $G$, although we are looking at recurrence
in $H$, the patterns we are looking at come from a F\o{}lner sequence
in $G$, which is not contained in $H$. We emphasize that these results
do not follow from an application of our recurrence theorems to the
$H$-process which arises from the restriction of the $G$-action
to $H$; in fact the $H$-process derived from an ergodic $G$-process
need not be ergodic, and even if it is it need not have the same entropy
as the original process. This observation means that we cannot weaken
the requirement that $F_{n}$ be a F\o{}lner sequence. If for example
$\{F_{n}\}$ were a F\o{}lner sequence in some proper subgroup $H$
of $G$, our results as applied to the $H$-process and the ergodic
component to which $x$ belongs would show that $U_{*},U^{*}$ etc.
do converge, but their value depends on $x$ and may be different
from the entropy of the original process.

\subsection{\label{subsec:bad-window-sequences}Return times and densities may
misbehave for {}``nice'' window sequences which aren't filling}

It is natural to wonder which of the theorems about return times and
densities are true under weaker conditions than those stated, and
in particular whether the requirement that $\left\{ W_{n}\right\} $
be filling may not be replaced with a weaker condition. One condition
which has been used successfully in the proof of other pointwise theorems
is that of being a Templeman sequence, ie an increasing F\o{}lner
sequence with $|W_{n}^{-1}W_{n}|\leq C|W_{n}|$ (note that this is
implies temperedness of the sequence $\left\{ W_{n}\right\} $). In
this section we show that even a stronger condition is not enough,
and that the almost-certain bound \[
\liminf _{k}\inf _{n}-\frac{1}{|F_{k}|}\log U_{k,n}^{(F,F)}(x)\leq h\]
 may fail even if $\left\{ F_{n}\right\} $ is a symmetric, increasing
Templeman sequence of all orders, ie for all $d$ there is a constant
$C_{d}$ such that $|(F_{n})^{d}|\leq C_{d}|F_{n}|$ for all $n$.

Let $S^{1}=\left\{ z\in \mathbb{C}\, :\, |z|=1\right\} $ be the unit
circle in the complex plane and let $Tz=e^{2\pi i\theta }z$ be an
irrational rotation, ie $\theta \in \mathbb{R}\setminus \mathbb{Q}$.
Let $X:S^{1}\rightarrow \left\{ 0,1\right\} $ be the map $z\mapsto \, \textrm{sgn}\, \textrm{Im}\, z$,
and consider the process $X_{n}=T^{n}X$ defined on $S^{1}$ when
$S^{1}$ is equipped with Lebesgue measure. It is well known that
this process is ergodic and has zero entropy. Our goal here is to
construct a sequence of increasing, symmetric finite sets $F_{n}\subseteq \mathbb{Z}$
such that $\left\{ F_{n}\right\} $ is Templeman of all orders and
$|F_{n}|=2n$, but for which there is a sequence of indices $k_{i}$
such that $T_{k_{i}}^{(F,F)}(z)>2^{|F_{k_{i}}|}$ for every $z\in S^{1}$,
so $\frac{1}{|F_{k}|}\log T_{k}^{(F,F)}$ does not converge even in
probability to the entropy. This implies that the asymptotic upper
bound for $\inf _{n}-\frac{1}{|F_{k}|}\log U_{k,n}^{(F,F)}$ fails
as well.

Let $\varepsilon _{n}\searrow 0$ be a decreasing sequence of positive
numbers such that for every $z\in S^{1}$, if $0<|k|<n$ then $|T^{k}z-z|>\varepsilon _{n}$.
For any $z_{0}\in S^{1}$ one can find such an $\varepsilon _{n}$,
since $T$ has no periodic orbits. Since $T$ is an isometry of $S^{1}$
such a choice of $\varepsilon _{n}$ is good for all $z\in S^{1}$
if it is good for $z_{0}$.

Let $\ell _{n}$ be an increasing sequence of positive integers such
that any realization $(X_{-\ell _{n}}(z),\ldots ,X_{\ell _{n}}(z))$
determines $z$ up to a distance of $\varepsilon _{n}$. More precisely,
we require that if $z,z'\in S^{1}$ satisfy $X_{k}(z)=X_{k}(z')$
for all $-n\leq k\leq n$ then $|z-z'|<\varepsilon _{n}$. Such $\ell _{n}$
exist because if we denote by $\mathcal{P}$ the partition of $S^{1}$
determined by $X$ then $\lor _{-\ell }^{\ell }T^{n}\mathcal{P}$
partitions $S^{1}$ into segments whose length tends to $0$ as $\ell \rightarrow 0$.

Finally, write \[
J_{n}=\left\{ k\in \mathbb{Z}\, :\, |z-T^{k}z|>\varepsilon _{n}\textrm{ }(\textrm{for any or all }z\in S^{1})\right\} .\]
 Note that $J_{n}$ are symmetric (because $T$ is an isometry) and
that since $\varepsilon _{n}$ are decreasing, $J_{n}\subseteq J_{n+1}$.
Also note that if $I=[a;b]\subseteq \mathbb{Z}$ then\[
k,m\in I\setminus J_{n}\, \Rightarrow \, |z-T^{k-m}z|=|T^{k}z-T^{m}z|<2\varepsilon _{n}\]
 and so if $\varepsilon _{n'}>2\varepsilon _{n}$ then $n'<|k-m|$.
Thus if $|I|<n'$ then $|J_{n}\cap I|\geq |I|-1$. It follows that
the relative density of $J_{n}$ in any segment $I$ of length $n'-1$
is at least $1-\frac{1}{n'-1}$.

We define the sets $F_{n}$ inductively. At each stage we obtain $F_{n}$
from $F_{n-1}$ by appending to $F_{n-1}$ some symmetric pair of
numbers $\pm m$. Given $F_{n}$ and a symmetric set $E\subseteq \mathbb{Z}$,
when we will say that we {}``add $E$ to $F_{n}$'' we will mean
that over several stages we add to $F_{i}$ the smallest pair of numbers
(in absolute value) in $E\setminus F_{i}$ until $E$ is exhausted
or some condition is met.

We start with $F_{0}=\emptyset $ and add to it the segment $[-\ell _{1};\ell _{1}]$
to obtain the sequence $F_{1},\ldots ,F_{k_{1}}$. We now add $J_{1}$
until we obtain a set $F_{n_{1}}$ such that $|F_{n_{1}}|>2^{|F_{k_{1}}|}$.
Now we add $[-\ell _{2};\ell _{2}]$ to $F_{n_{1}}$ and obtain $F_{k_{2}}$;
then add $J_{2}$ to $F_{k_{2}}$ until we obtain a set $F_{n_{2}}$
with $|F_{n_{2}}|>2^{|F_{k_{2}}|}$. We proceed in this manner, alternately
adding $[\ell _{-i};\ell _{i}]$ to obtain $F_{k_{i}}$ and then $J_{i}$
to obtain $F_{n_{i}}$ with $|F_{n_{i}}|>2^{|F_{k_{i}}|}$.

The sequence $\left\{ F_{n}\right\} $ thus constructed is clearly
increasing and symmetric, and satisfies $|F_{n}|=2n$. Next, we claim
that $T_{k_{i}}(z)>2^{|F_{k_{i}}|}$ for every $z\in S^{1}$. To see
that this is true for all $i$, note that $z(F_{k_{i}})$ determines
$z$ up to $\varepsilon _{i}$ because $[-\ell _{i};\ell _{i}]\subseteq F_{k_{i}}$,
while for every $m\in F_{n_{i}}\setminus F_{k_{i}}$ we have $m\in J_{i}$,
so $|z-T^{m}z|>\varepsilon _{i}$. Therefore $T^{m}z([-\ell _{i};\ell _{i}])\neq z([-\ell _{i};\ell _{i}])$
and this forces $T^{m}z(F_{k_{i}})\neq z(F_{k_{i}})$. Since $|F_{n_{i}}|>2^{|F_{k_{i}}|}$,
this proves the claim.

As for $\left\{ F_{n}\right\} $ being Templeman, it suffices to show
that $\lim _{n\rightarrow \infty }\frac{|F_{n}|}{\max F_{n}-\min F_{n}}=1$;
this implies both that $\left\{ F_{n}\right\} $ is a F\o{}lner sequence,
and that $|F_{n}^{d}|\leq C_{d}|F_{n}|$ for some constant $C_{d}$,
for all $n$. But this follows from our remarks about the relative
density of $J_{n}$ in certain segments.

We remark that it seems to be more difficult to construct an example
of such a sequence $\left\{ F_{n}\right\} $ for which the lower bound
for $U_{k,n}$ fails. Perhaps the lower bound $U_{*}\geq h$ holds
more generally than the upper bound, as happens in the case of $T_{k}$.
we have so far been unable to determine whether this is so.

\section{\label{sec:upper-bounds}Upper bounds for return times and densities}

In this section we prove the upper bounds $T^{*}\leq h$ and $U^{*}\leq h$
as described in the introduction.

\subsection{\label{subsec:upper-bound-for-U} The upper bound for $U^{*}$}

We begin with some further remarks about the definition of $U_{k,n}$.
Recall that in section \ref{subsec:recurrence-densities} we defined
$U_{k,n}$ to be size of the maximal collection of repetitions of
$x(F_{k})$ in $x(W_{n})$ satisfying (a) the collection includes
the trivial repetition $x(F_{k})$ itself, and (b) the members of
the collection are pairwise disjoint. If we drop these restriction
we obtain the quantity \[
V_{k,n}(x)=V_{k,n}^{(F,W)}(x)=\frac{1}{|W_{n}|}\#\{f\in W_{n}\, :\, x(F_{k})\textrm{ repeats at }f\}\]
 which is related to $U_{k,n}$ by \begin{equation}
U_{k,n}\leq V_{k,n}\leq |F_{k}^{-1}F_{k}|U_{k,n}\leq |F_{k}|^{2}U_{k,n}.\label{eq:U-and-U-bar-equivalence}\end{equation}
 The middle inequality is due to the fact that from any collection
of right translates of $F_{k}$ we can obtain a pairwise disjoint
sub-collection at least $\frac{1}{|F_{k}^{-1}F_{k}|}$ the size of
the original collection which includes any single prescribed member
of the original collection (which in our case we take to be the set
$F_{k}$, representing the original pattern). From this relation we
see that in order to prove that $U^{*}\leq h$, it suffices to prove
that $V^{*}=\limsup _{n}\, \sup _{k}-\frac{1}{|F_{k}|}\log V_{k,n}\leq h$.

The key lemma which makes filling sequences special is the following,
which states that if a set $E$ has a low {}``local'' density with
respect to $\left\{ W_{n}\right\} $ then it has a low {}``global''
density:

\begin{lem}
\label{lemma:basic-filling-lemma}Suppose $\left\{ W_{n}\right\} $
is a $C$-filling sequence and $E\subseteq G$ finite. If for each
$f\in E$ there is an index $n(f)$ such that $\frac{1}{|W_{n(f)}|}|W_{n(f)}f\cap E|\leq \alpha $
then \[
\frac{|E|}{|\bigcup _{f\in E}W_{n(f)}f|}\leq \frac{\alpha }{C}.\]

\end{lem}
Note that we always have $0<C\leq 1$.

\begin{proof}
Order $E$ by decreasing value of $n(\cdot )$: $E=\left\{ f_{1},\ldots ,f_{I}\right\} $
and write $n(f_{i})=n(i)$. We select a subsequence $\left\{ f_{i(j)}\right\} \subseteq \left\{ f_{1},\ldots ,f_{I}\right\} $
inductively. Set $i(1)=1$ and let $i(j)$ be the first index such
that $f_{i(j)}\notin \bigcup _{k<j}W_{n(k)}f_{k}$; the process ends
after $J$ steps when $E\subseteq \bigcup _{j=1}^{J}W_{n(i(j))}f_{i(j)}$
(recall that $1_{G}\in W_{n}$ for each $n$, since $\left\{ W_{n}\right\} $
is filling). Clearly $\left\{ W_{n(i(j))}f_{i(j)}\right\} _{j=1}^{J}$
is incremental. Now \[
|E|\leq \sum _{j=1}^{J}|W_{n(i(j))}f_{i(j)}\cap E|\leq \alpha \sum _{j=1}^{J}|W_{n(i(j))}|\leq \frac{\alpha }{C}|\bigcup _{j=1}^{J}W_{n(i(j))}f_{i(j)}|\]
 (the middle inequality is because by hypothesis $|W_{n(i)}f_{i}\cap E|\leq \alpha |W_{n(i)}|$
for every $i=1,\ldots ,I$ and the last because $\left\{ W_{n}\right\} $
is $C$-filling). Since \[
\bigcup _{j=1}^{J}W_{n(i(j))}f_{i(j)}\subseteq \bigcup _{i=1}^{I}W_{n(i)}f_{i}\]
 the lemma follows. 
\end{proof}
We apply lemma \ref{lemma:basic-filling-lemma} in the following situation:
Consider an ergodic $G$-process defined on $(\Omega ,\mathcal{F},P)$
with entropy $h$. Suppose $x\in \Omega $, $k<L<<N$ are fixed, and
$\varphi :F_{k}\rightarrow \Sigma $ is some $F_{k}$ pattern. Let
$E$ be the set of repetitions of $\varphi $ in $x(F_{N})$ with
the additional property that near them $\varphi $ repeats not frequently
enough with respect to $V_{k,n}$: \[
E=\left\{ f\in F_{N}\, \left|\begin{array}{c}
 fx\in [\varphi ]\textrm{ and there exists }n(f)\leq L\\
 \textrm{s}.\textrm{t}.\textrm{ }V_{k,n(f)}(fx)<2^{-(h+\varepsilon )|F_{k}|}\end{array}
\right.\right\} .\]
 Since for any $f\in E$ we have $fx\in [\varphi ]$ and $\frac{1}{|W_{n(f)}|}|E\cap W_{n(f)}f|<2^{-(h+\varepsilon )|F_{k}|}$,
from the lemma we see that \[
\frac{1}{|\cup _{f\in E}W_{n(f)}f|}|E|<\frac{1}{C}\cdot 2^{-(h+\varepsilon )|F_{k}|}.\]
 If we assume that $N$ is large enough so that $|(\cup _{n\leq L}W_{n})F_{N}|\leq 2|F_{N}|$
(this happens eventually since $\left\{ F_{n}\right\} $ is a F\o{}lner
sequence), then since $\cup _{f\in E}W_{n(f)}f\subseteq (\cup _{n\leq L}W_{n})F_{N}$,
we see that \[
\frac{1}{|F_{N}|}|E|\leq \frac{2}{|\cup _{f\in E}W_{n(f)}f|}|E|<\frac{2}{C}\cdot 2^{-(h+\varepsilon )|F_{k}|}.\]

Given some $\Phi \subseteq \Sigma ^{F_{k}}$ and repeating this argument
for every $\varphi \in \Phi $ we have\[
\frac{1}{|F_{N}|}\#\left\{ f\in F_{N}\, \left|\begin{array}{c}
 fx\in [\Phi ]\textrm{ and there exists }n(f)\leq L\\
 \textrm{s}.\textrm{t}.\textrm{ }V_{k,n(f)}(fx)<2^{-(h+\varepsilon )|F_{k}|}\end{array}
\right.\right\} <\frac{2}{C}\cdot 2^{-(h+\varepsilon )|F_{k}|}\cdot |\Phi |.\]

\begin{thm}
\label{thm:upper-bound-for-U}If $\{F_{n}\}$ is a tempered F\o{}lner
sequence and if $\{W_{n}\}$ is filling, then \[
U^{*}(x)=\limsup _{k\rightarrow \infty }\, \sup _{n}-\frac{1}{|F_{k}|}\log U_{k,n}^{(F,W)}(x)\leq h\]
 almost surely. 
\end{thm}
\begin{proof}
It is enough to prove that \[
V^{*}=\limsup _{k}\, \sup _{n}-\frac{1}{|F_{k}|}\log |V_{k,n}|\leq h\]
 almost surely. Suppose to the contrary that for some $\varepsilon >0$
and measurable non-null set $B\subseteq \Omega $, for every $x\in B$
\[
\limsup _{k\rightarrow \infty }\, \sup _{n}-\frac{1}{|F_{k}|}\log V_{k,n}^{(F,W)}(x)>h+\varepsilon .\]
 Let $P(B)>p>0$. Then from the SMB theorem, there is a measurable
subset $B_{0}\subseteq B$, $P(B_{0})>p$, a sequence $\left\{ \Phi _{k}\right\} $
of sets $\Phi _{k}\subseteq \Sigma ^{F_{k}}$ and integer $L_{1}$
such that for $k>L_{1}$ it holds that $|\Phi _{k}|\leq 2^{(h+\varepsilon /2)|F_{k}|}$,
and $x(F_{k})\in \Phi _{k}$ for any $x\in B_{0}$. It may further
be assumed that for some $L_{2}>>L_{1}$, for any $x\in B_{0}$ there
are indices $k(x),n(x)$ with $L_{1}\leq k(x)\leq L_{2}$ and $n(x)\leq L_{2}$
such that $V_{k(x),n(x)}(x)<2^{-(h+\varepsilon )|F_{k(x)}|}$.

Choose $N$ very large with respect to $L_{2}$ so that $F_{N}$ is
very $W_{L_{2}}$-invariant, and large enough so that there is an
$x\in \Omega $ such that\[
\frac{1}{|F_{N}|}\#\{f\in F_{N}\, :\, fx\in B_{0}\}>p\]
 (by the ergodic theorem such $x$'s are bound to exist for $N$ large
enough).

Fix $x\in \Omega $ for which the above holds. For $f\in F_{N}$ such
that $fx\in B_{0}$, write $k(f)=k(fx)$, $n(f)=n(fx)$. Define\[
E=\{f\in F_{N}\, :\, fx\in B_{0}\},\]
 \[
E_{k}=\{f\in F_{N}\, :\, fx\in B_{0}\textrm{ and }k(f)=k\}.\]
 Clearly $E\subseteq (\cup _{k=L_{1}}^{L_{2}}E_{k})$, so by our choice
of $x$ we should have\[
\frac{1}{|F_{N}|}|\cup E_{k}|\geq p.\]

On the other hand, according to the discussion after lemma \ref{lemma:basic-filling-lemma},
\[
\frac{1}{|F_{N}|}|E_{k}|\leq \frac{2}{C}2^{-(h+\varepsilon )|F_{k}|}|\Phi _{k}|\leq \frac{2}{C}(2^{-\varepsilon /2})^{|F_{k}|}\]
 so that \begin{equation}
\frac{1}{|F_{N}|}|\cup E_{k}|=\frac{1}{|F_{N}|}\sum _{L_{1}<k<L_{2}}|E_{k}|\leq C'\cdot (2^{-\varepsilon /2})^{F_{L_{1}}}\label{eq:U-lower-bound-punchline}\end{equation}
 for some constant $C'$ depending only on $\varepsilon $, so if
$L_{1}$ was chosen large enough with respect to $\varepsilon $ (recall
our assumption that that $|F_{n}|\geq n$), this is less than $p$,
a contradiction. 
\end{proof}
Esssencially the same proof gives us

\begin{thm}
\label{thm:upper-bound-for-U-weak-version}Let $\{F_{n}\}$ be a tempered
F\o{}lner sequence in $G$. If $\left\{ W_{n}\right\} $ is a sequence
such that $W_{n}$ is $C_{n}$-filling and $\left\{ n(k)\right\} $
is a sequence of integers such that $\sum _{k=1}^{\infty }\frac{2^{-\alpha |F_{k}|}}{C_{n(k)}}<\infty $
for every $\alpha >0$, then for almost every $x$, \[
\limsup _{k\rightarrow \infty }\, -\frac{1}{|F_{k}|}\log U_{k,n(k)}^{(F,W)}(x)\leq h.\]

\end{thm}
The proof uses the fact that lemma \ref{lemma:basic-filling-lemma}
has an analogue when $\left\{ W_{n}\right\} $ is replaced with a
single $C$-filling set $W$, namely that if for a finite subset $E\subseteq G$
we have $\frac{|Wf\cap E|}{|W|}<\alpha $ for every $f\in E$ then
$\frac{|E|}{|WE|}\leq \frac{\alpha }{C}$ (to see this, apply the
lemma with $W_{n}=W$ ). Using this the proof now follows exactly
the lines of the proof of theorem \ref{thm:upper-bound-for-U}; the
relation $\sum \frac{2^{-\alpha |F_{k}|}}{C_{n(k)}}<\infty $ is used
to show that the middle term of equation (\ref{eq:U-lower-bound-punchline})
is small when $L_{1}$ is large enough.

\subsection{\label{subsec:upper-bound-for-T}The upper bound for $T^{*}$}

In order for $T_{k}^{(F,W)}$ to be meaningful the window sequence
$\left\{ W_{n}\right\} $ must grow slowly enough to detect a too-soon
return. What slowly means here is that for every $h$ there exist
window sets of size approximately $2^{h|F_{k}|}$ for large enough
$k$. One condition which ensures this without reference to the sequence
$\left\{ F_{n}\right\} $ is that $\frac{\log |W_{n+1}|}{\log |W_{n}|}\rightarrow 1$.
With this condition in place, one can deduce the bound for $T^{*}$
from the bound for $U^{*}$:

\begin{lem}
\label{lemma:U-upper-bound-gives-T-upper-bound}Suppose $\left\{ W_{n}\right\} $
is an increasing sequence such that $\frac{\log |W_{n+1}|}{\log |W_{n}|}\rightarrow 1$.
For $x\in \Omega $, if $U^{*}(x)\leq h$ then $T^{*}(x)\leq h$. 
\end{lem}
\begin{proof}
If $T_{k}(x)=|W_{R_{k}(x)}|>2^{(h+\varepsilon )|F_{k}|}$ for infinitely
many $k$, then the growth condition implies that $|W_{R_{k}(x)-1}|\geq 2^{(h+\varepsilon /2)|F_{k}|}$
for infinitely many $k$. But for such a $k$ the only occurence of
$x(F_{k})$ in $W_{R_{k}(x)-1}$ is the original pattern, so $|U_{k,R_{k}(x)-1}(x)|\leq 2^{-(h+\varepsilon /2)|F_{k}|}$
for infinitely many $k$, which is impossible because $U^{*}(x)\leq h$. 
\end{proof}
\begin{thm}
\label{thm:upper-bound-for-T-strong-version}Let $\left\{ W_{n}\right\} $
be a quasi-filling sequence satisfying $\frac{\log |W_{n+1}|}{\log |W_{n}|}\rightarrow 1$.
Then for any tempered F\o{}lner sequence $\left\{ F_{n}\right\} $,
$T_{(F,W)}^{*}\leq h$ almost surely. 
\end{thm}
\begin{proof}
Let $\varepsilon >0$. Define \[
n(k)=\min \left\{ n:|W_{n}|>2^{(h+\varepsilon )|F_{k}|}\right\} .\]
For $\alpha >0$ we have, since $\{W_{n}\}$ is quasi-filling, that
\[
\sum \frac{2^{-\alpha |F_{k}|}}{C_{n(k)}}\leq \sum (\frac{1}{C_{n(k)}}2^{(h+\varepsilon )|F_{k}|})^{-\alpha /(h+\varepsilon )}\leq \sum \frac{|W_{n(k)}|^{-\alpha /(h+\varepsilon )}}{C_{n(k)}}<\infty \]
 and this now gives, via theorem \ref{thm:upper-bound-for-U-weak-version},
that\[
\limsup _{k}\frac{1}{|F_{k}|}\log U_{k,n(k)}^{(F,W)}(x)\leq h\]
 almost surely. This in turn implies that $T^{*}\leq h$ almost surely,
by reasoning like that of the lemma \ref{lemma:U-upper-bound-gives-T-upper-bound}
above. 
\end{proof}

\section{\label{sec:lower-bounds}Lower bounds on return times and densities}

In this section we prove the bounds $T_{*}\geq h$ and $U_{*}\geq h$.
The bound $T_{*}\geq h$ is proved by constructing a code for the
process which beats entropy if $T_{*}\not \geq h$. This is close
in spirit to the proof of $T_{*}\geq h$ originally given by Ornstein
and Weiss for the cases $G=\mathbb{Z}$ and $G=\mathbb{Z}^{d}$, but
some new combinatorial machinery is needed to make the construction
possible for general groups and general window sets. We first give
a brief description of the coding ideas we will need, and then explain
the combinatorics and detailed construction of the code. We conclude
with the proof of the bound $U_{*}\geq h$, which relies both on the
filling ideas used in the last section and on the lower bound $T_{*}\geq h$.

\subsection{\label{subsec:coding} Coding and entropy}

For completeness we take a detour and discuss the connection between
entropy and efficient coding of a process.

A useful characterization of the entropy of a process is as the lower
bound of all achievable coding rates for the process. An $F_{n}$-\emph{code}
is a map $c:\Sigma ^{F_{n}}\rightarrow \{0,1\}^{*}=\cup _{n\geq 0}\{0,1\}^{n}$.
Thus $c$ encodes each pattern $\varphi :F_{n}\rightarrow \Sigma $
in a binary string (=sequence) $c(\varphi )$, which is called the
\emph{codeword} associated with $\varphi $. If $c$ is an injection
the coding is said to be \emph{invertible or} \emph{faithful}. For
$x\in \Omega $ and an $F_{n}$-code $c$, we write $c(x)$ instead
of $c(x(F_{n}))$. If $c_{n}$ is an $F_{n}$-code, we say that $\left\{ c_{n}\right\} $
is an $\left\{ F_{n}\right\} $-code. $\left\{ c_{n}\right\} $ is
invertsble if every $c_{n}$ is.

The number of bits (=symbols) in the codeword $c(\varphi )$ of $\varphi $
is called the \emph{length} of the codeword $c(\varphi )$ and is
denoted by $\ell (c(\varphi ))$. Note that the codeword length of
different patterns may vary. For an $F_{n}$-code $c$ and process
$\left\{ X_{g}\right\} $, the average number of \emph{bits per symbol}
used by the code is the quantity $\rho (c)=\frac{1}{|F_{n}|}\int \ell (c(x))dP(x)$;
this is the average number of bits used in the $c$-encoding of a
random pattern $x(F_{n})$, normalized by the number of symbols in
the $F_{n}$-pattern $x(F_{n})$.

For $r\geq 0$, we say that $r$ is an \emph{achievable coding rate}
for the process \emph{$\left\{ X_{g}\right\} $} if there exists an
invertable $\left\{ F_{n}\right\} $-code $\left\{ c_{n}\right\} $
such that $\limsup \rho (c_{n})\leq r$. This means we can represent
arbitrarily large patterns from the process using approximately $r$
bits for every symbol of the pattern. The \emph{coding rate} of the
process $\left\{ X_{g}\right\} $ is the infimum of achievable coding
rates: \[
R=\inf \left\{ r\geq 0\, :\, r\textrm{ is an achievable coding rate for the process }\left\{ X_{g}\right\} \right\} .\]

It is not hard to see that the coding rate of a process is the same
as the process entropy $h$ (incidentally, this shows that $R$ does
not depend on the sequence $\left\{ F_{n}\right\} $). To see that
$R\geq h$, let $\{c_{n}\}$ be a sequence of invertable codes such
that $\limsup _{n\rightarrow \infty }\rho (c_{n})<r$. Then there
is some $\varepsilon >0$ such that for infinitely many $n$ there
is a measurable set $A_{n}\subseteq \Omega $ with $P(A_{n})\geq \varepsilon $,
and $\frac{1}{|F_{n}|}\ell (c_{n}(x))<r$ for every $x\in A_{n}$.
Since $c_{n}$ is an injection, and for every $x\in A_{n}$ we have
$\ell (c_{n}(x))\leq r|F_{n}|$, we have \[
\{c_{n}(x)\, :\, x\in A_{n}\}\subseteq \bigcup _{k<r|F_{n}|}\left\{ 0,1\right\} ^{k}\]
 so every $x\in A_{n}$ belongs to at most one $\sum _{k\leq r|F_{n}|}2^{k}\leq 2\cdot 2^{r|F_{n}|}$
$F_{n}$-atoms. The Shannon-McMillan theorem now implies that $h\leq r$.

To see that $h\geq R$ we formulate the following simple lemma, which
enables us to turn a non-faithful code sequence into a faithful code
sequence.

\begin{lem}
\label{lemma:making-faithful-codes-from-nonfaithful-codes}Suppose
that $\{c_{n}\}$ is a sequence of $F_{n}$-codes and $\Phi _{n}\subseteq \Sigma ^{F_{n}}$
is such that $c_{n}|_{\Phi _{n}}$ is an injection and $\frac{1}{|F_{n}|}\ell (c_{n}(\varphi ))\leq r$
for every $\varphi \in \Phi _{n}$. Assume also that $P([\Phi _{n}])\rightarrow 1$.
Then $R\leq r$. 
\end{lem}
\begin{proof}
Let $e_{n}$ be an enumeration coding of $\Sigma ^{F_{n}}$, ie it
assigns to $\varphi \in \Sigma ^{F_{n}}$ a binary string representing
the index of $\varphi $ in some fixed enumeration of $\Sigma ^{F_{n}}$.
The code $e_{n}$ is faithful. Now define a code $\{\hat{c}_{n}\}$
by\[
\hat{c}_{n}(\varphi )=\left\{ \begin{array}{cc}
 0c_{n}(\varphi ) & \varphi \in \Phi _{n}\\
 1e_{n}(\varphi ) & \varphi \notin \Phi _{n}\end{array}
\right.\]
 $\{\hat{c}_{n}\}$ is clearly an invertible code. Moreover, \[
\frac{1}{|F_{n}|}\ell (\hat{c}_{n}(\varphi ))=\left\{ \begin{array}{cc}
 \frac{1}{|F_{n}|}+\frac{1}{|F_{n}|}\ell (c_{n}(\varphi )) & \varphi \in \Phi _{n}\\
 \frac{1}{|F_{n}|}+\left\lceil \log |\Sigma |\right\rceil  & \varphi \notin \Phi _{n}\end{array}
\right.\]
 so, \begin{eqnarray*}
\rho (\widehat{c}_{n}) & = & \frac{1}{|F_{n}|}\int \ell (\hat{c}_{n}(x))dP\\
 & \leq  & \frac{1}{|F_{n}|}+r\cdot P([\Phi _{n}])+\left\lceil \log |\Sigma |\right\rceil P(\Omega \setminus [\Phi _{n}])\\
 & \rightarrow  & r
\end{eqnarray*}
as desired.
\end{proof}
If we choose $\Phi _{n}\subseteq \Sigma ^{F_{n}}$ in such a way that
$\frac{1}{|F_{n}|}\log |\Phi _{n}|\rightarrow h$ and $P([\Phi _{n}])\rightarrow 1$,
choose $c_{n}$ to be an enumeration coding of $\Phi _{n}$ and define
it arbitrarily of $\Sigma ^{F_{n}}\setminus \Phi _{n}$, then the
lemma implies that $\rho (\widehat{c}_{n})\leq h+\varepsilon $ for
every $\varepsilon >0$, so $R\leq h$.

We conclude this section with a combinatorial lemma which will be
used later. Its importance is that it enables us to describe small
random subsets of $F_{n}$ at a low rate by using an enumeration code.

\begin{lem}
\label{lemma:number-of-small-subsets} Let $0<\lambda <1/2$ and $F$
a finite set. The number of subsets of $F$ whose size is at most
$\lambda |F|$ is bounded by $2^{\delta (\lambda )|F|}$, where $\lim _{\lambda \rightarrow 0}\delta (\lambda )=0$. 
\end{lem}
\begin{proof}
The number of such subsets is \[
\sum _{m<\lambda |F_{n}|}{{|F_{n}| \choose m}}\leq 2^{(-\lambda \log \lambda -(1-\lambda )\log (1-\lambda ))|F|}.\]
 The bound follows from Stirling's formula; or see \cite[p. 52]{Sh96}
for an elementary proof. 
\end{proof}

\subsection{\label{subsec:lower-bound-for-T-part-1}The lower bound for $T_{*}$
(Part I)}

In this section we begin the proof of

\begin{thm}
\label{thm:lower-bound-for-T} If $\{F_{n}\}$ is a tempered sequence
and $\{W_{n}\}$ is an increasing sequence of finite subsets of $G$,
then \[
T_{*}^{(F,W)}=\liminf _{k\rightarrow \infty }\frac{1}{|F_{k}|}\log T_{k}^{(F,W)}\geq h\]
 almost surely. 
\end{thm}
We will need some notation.

\begin{defn}
\label{notation:cover}Let $\{F_{n}\}$ be a fixed sequence. A \emph{cover
$\nu $} is a collection of sets of the form $F_{n}f$ for some $n\in \mathbb{N}$
and $f\in G$. The set $E\subseteq G$ of $f$'s such that some $F_{n}f$
is represented in $\nu $ is called the \emph{set of centers of $\nu $}
and is denoted by $\dom \nu $. We also say then that $\nu $ is a
cover \emph{over} E. 
\end{defn}
It would be more precise to say that a cover $\nu $ over $E$ is
a collection of pairs $(n,f)\in \mathbb{N}\times E$, because we wish
to allow more than one set $F_{n}f$ over each $f\in E$, and also
since the set $F_{n}f$ does not, in general, determine $n$ and $f$
uniquely. However, it is more convenient to think of $\nu $ as a
(multi)set with elements of the form $F_{n}f$, and this is what we
will do.

We write $\cup \nu $ for the union of members of $\nu $. When we
wish to specify the elements of $\nu $ by name, we will write $\cup _{\nu }F_{n}f$,
$\sum _{\nu }|F_{n}f|$, etc. This is imprecise because the index
$n$ depends on $f$ but we omit this when convenient.

Fix $\left\{ F_{n}\right\} $ and $\left\{ W_{n}\right\} $. Let $x\in \Omega $,
$f\in G$ and $n\in \mathbb{N}$, and suppose $R_{n}(fx)<\infty $.
There then exists an element $f'\in G$ at which $fx(F_{n})$ repeats,
and is closest to $f$ with respect to $\left\{ W_{n}\right\} $.
To be precise, there is an $f'$ such that

\begin{enumerate}
\item $f'x\in [fx(F_{n})]$ (ie the $F_{n}$-pattern at $f$ in $x$ repeats
in $f'$).
\item $f'f^{-1}\in W_{R_{n}(fx)}$ (ie the offset from $f$ to $f'$ is
in $W_{R_{n}(f(x))}$).
\end{enumerate}
In general there may be more than one $f'$ with these properties.
For fixed $x$ and $n$, we assume some choice has been made, and
denote this element by $f^{x,n}$.

For a set $F_{n}f$ such that $R_{n}(fx)<\infty $, we write $(F_{n}f)^{x}=F_{n}f^{n,x}$.
The (partial) map $F_{n}f\mapsto (F_{n}f)^{x}$ is called the \emph{displacement
arising from $x$.} Strictly speaking, the map $F_{n}f\mapsto (F_{n}f)^{x}$
is not well defined since the set $F_{n}f$ does not uniquely determine
$n$ and $f$. However, $n$ and $f$ will always be clear from the
context so no confusion should arise. Finally, for a cover $\nu $
such that $R_{n}(fx)<\infty $ for every $F_{n}f\in \nu $ we write
$\nu ^{x}=\left\{ (F_{n}f)^{x}\, :\, F_{n}f\in \nu \right\} $. Thus
$\nu ^{x}$ is a cover over the set $\left\{ f^{x,n}\, :\, F_{n}f\in \nu \right\} $.
One should think of the displacement $\nu ^{x}$ of $\nu $ as the
cover obtained from $\nu $ by moving each set $F_{n}f\in \nu $ to
one of the nearest repetitions of the pattern induced on it by $x$.

We will only be working with covers $\nu $ and $x\in \Omega $ such
that $R_{n}(fx)<\infty $ for every $F_{n}f\in \nu $, and avoid further
mention of this fact. 

We can now state the idea behind the proof of theorem \ref{thm:lower-bound-for-T}.
Suppose by way of contradiction that $P(T_{*}<h-\varepsilon )>p>0$.
We will try to construct an $F_{N}$-code $c_{N}$ which beats entropy.
By the ergodic theorem, for almost every $x$ and large $N$ at least
a $p$-fraction of the $f\in F_{N}$ are such that $T_{n}(fx)\leq 2^{(h-\varepsilon )|F_{n}|}$
for some $n=n(f)\in \mathbb{N}$. In order to describe the pattern
$x(F_{n}f)$ it is enough to describe the pattern $x(F_{n}f^{x,n})$
and the element $f^{x,n}f^{-1}$. Since the latter is in $W_{R_{n}(fx)}$
and this set is of size at most $2^{(h-\varepsilon )|F_{n}|}$, we
can describe it in approximately $(h-\varepsilon )|F_{n}|$ bits.
Thus if the pattern $x(F_{n}f^{x,n})$ were known, the pattern $x(F_{n}f)$
could be coded in $h-\varepsilon $ bits per symbol.

Let $E\subseteq F_{N}$ be the set of $f\in F_{N}$ for which there
exists such an $n(f)$ with $R_{n(f)}(fx)\leq h-\varepsilon $. Let
$\nu =\left\{ F_{n(f)}f\right\} _{f\in E}$. If it should happen that
(a) the collection $\nu $ is pairwise disjoint, (b) $\cup \nu $
is disjoint from $\cup \nu ^{x}$, (c) we can code $x(F_{N}\setminus \cup \nu )$
at, say $h+p\varepsilon /2$ bits per symbol, and (d) $\cup \nu $
is at least a $p$-fraction of $F_{N}$, then the discussion above
shows that we can code $x(\cup \nu )$ at around $h-\varepsilon $
bits per symbol. Given that the above hold, we can then code $x(F_{N})$
at a rate of $h-p\varepsilon /2$ bits per symbol.

Conditions (c) and (d) are not to difficult: (d) would follow from
(a) and the fact that $E$ is a $p$-fraction of $F_{N}$, and (c)
can be achieved because $F_{N}\setminus \cup \nu $ can be mostly
covered by entropy-typical sets when $N$ is large, and a standard
argument then gives that $x(F_{N}\setminus \cup \nu )$ can be coded
at a little more than $h$ bits per symbol. The difficulty in realizing
this outline is that in (a) and (b) cannot in general be achieved.

The proof of the return times theorem for $\mathbb{Z}$ and $\mathbb{Z}^{d}$
given in \cite{OW02} avoids this problem by allowing the elements
of $\nu $ and $\nu ^{x}$ to intersect arbitrarily, and employs a
more elaborate coding scheme which uses the order structure of $\mathbb{Z}^{d}$
in an essential way. This method doesn't generalize to other groups.

For our proof we will drop the disjointness requirements (a) and (b)
and instead require only {}``almost disjointness''. To be precise,

\begin{defn}
\label{def:epsilon-disjointness} Let $\varepsilon >0$. A sequence
$\{H_{i}\}_{i\in I}$ of finite subsets of $G$ are said to be \emph{$\varepsilon $-disjoint}
if there exist subsets $H'_{i}\subseteq H_{i}$ such that $\{H'_{i}\}_{i\in I}$
are pairwise disjoint and $|H'_{i}|\geq (1-\varepsilon )|H_{i}|$.
We say that a cover $\nu $ is $\varepsilon $-disjoint if the collection
$\{F_{n}f\in \nu \}$ is $\varepsilon $-disjoint. 
\end{defn}
If $\{H_{i}\}_{i=1}^{I}$ has the property that for each $i\leq I$,
$|H_{i}\cap (\cup _{j<i}H_{j})|\leq \varepsilon |H_{i}|$ then the
collection is $\varepsilon $-disjoint; simply set $H'_{i}=H_{i}\setminus \cup _{j<i}H_{j}$.

An important property of $\varepsilon $-disjoint collections is that
the size of their union is almost the sum of their sizes. To be precise,
let $\{H_{i}\}_{i\in I}$ be $\varepsilon $-disjoint and $H'_{i}\subseteq H_{i}$
as in the definition. Then we have \[
|\bigcup _{i\in I}H_{i}|\geq |\bigcup _{i\in I}H'_{i}|=\sum _{i\in I}|H'_{i}|\geq (1-\varepsilon )\sum _{i\in I}|H_{i}|.\]

Clearly for two $\varepsilon $-disjoint collections $\{D_{i}\}$
and $\{E_{j}\}$, if $D_{i}\cap E_{j}=\emptyset $ for every $i,j$
then $\{D_{i}\}\cup \{E_{j}\}$ is an $\varepsilon $-disjoint collection;
and so for any funion of pairwise disjoint $\varepsilon $-disjoint
collections.

\begin{defn}
\label{def:epsilon-disjointness-relative-to-displacement}For $x\in \Omega $
and a cover $\nu =\left\{ F_{n(i)}f_{i}\right\} _{i=1}^{I}$ over
$E=\left\{ f_{1},\ldots ,f_{I}\right\} $, we say that the pair $(\nu ,\nu ^{x})$
are $\varepsilon $-disjoint if for every $1\leq i\leq I$, \[
\left|F_{n(i)}f_{i}\cap \bigcup _{j<i}\left(F_{n(j)}f_{j}\cup \left(F_{n(j)}f_{j}\right)^{x}\right)\right|\leq \varepsilon |F_{n(i)}f_{i}|.\]

\end{defn}
According to the discussion before the definition, if $(\nu ,\nu ^{x})$
is $\varepsilon $-disjoint then $\nu $ is $\varepsilon $-disjoint.
Note that in order for the statement that $(\nu ,\nu ^{x})$ is $\varepsilon $-disjoint
to be meaningful, the set of centers of $\nu $ must be linearly ordered.

The point of this definitions is the following:

\begin{lem}
\label{lemma:pair-decoding}Let $0<\varepsilon <\frac{1}{2}$, $x\in \Omega $
and suppose $\nu =\left\{ F_{n(i)}f_{i}\right\} _{i=1}^{I}$ is a
cover over $E=\{f_{1},\ldots ,f_{I}\}$ such that $(\nu ,\nu ^{x})$
is $\varepsilon $-disjoint. Then there exists a set $A\subseteq \cup \nu $
satisfying $|A|<2\varepsilon |\cup \nu |$ and such that, given the
patterns $x(\cup \nu ^{x}\, \setminus \, \cup \nu )$ and $x(A)$,
and given the map $f_{i}\mapsto f^{x,n(i)}$, we can reconstruct the
pattern $x(\cup \nu )$. 
\end{lem}
\begin{proof}
Set \[
H_{i}=F_{n(i)}f_{i}\setminus \bigcup _{j<i}(F_{n(j)}f_{j})^{x}.\]
 Since $(\nu ,\nu ^{x})$ is $\varepsilon $-disjoint, $|H_{i}|\geq (1-\varepsilon )|F_{n(i)}|$,
and for $j>i$ we have $(F_{n(i)}f_{i})^{x}\cap H_{j}=\emptyset $.
Define \[
A=(\cup \nu )\setminus \bigcup _{i=1}^{I}H_{i}.\]
 Using the fact that $\nu $ is $\varepsilon $-disjoint (since $(\nu ,\nu ^{x})$
is) , \[
|A|<\varepsilon \sum |F_{n(i)}f_{i}|\leq \frac{\varepsilon }{1-\varepsilon }|\cup \nu |\leq 2\varepsilon |\cup \nu |.\]
 It remains to show that given $x(A)$, $x(\cup \nu ^{x}\, \setminus \, \cup \nu )$
and the map $f\mapsto f^{x,n(f)}$ we can deduce $x(\cup \nu )$.

We claim that for any $i$ if we know $x(\cup _{j<i}F_{n(j)}f_{j})$
then we can deduce $x(F_{n(i)}f_{i})$; this will complete the proof.
Since $f_{i}^{x,n(i)}x(F_{n(i)})=f_{i}x(F_{n(i)})$, and we know $f_{i}$
and $f_{i}^{x,n(i)}$, it suffices to show that we can deduce $x((F_{n(i)}f_{i})^{x})$.
Now since by assumption we already know $x(\cup _{j<i}F_{n(j)}f_{j})$,
and since $F_{n(i)}f_{i}\cap (F_{n(i)}f_{i})^{x}=\emptyset $, all
we need do is find $x((F_{n(i)}f_{i})^{x}\cap \bigcup _{j>i}F_{n(j)}f_{j})$.
But this is known, since $(F_{n(i)}f_{i})^{x}\cap \bigcup _{j>i}F_{n(j)}f_{j}\subseteq A$,
and we know $x(A)$. 
\end{proof}
In order to use the lemma \ref{lemma:pair-decoding} to produce good
codes we need to know how to produce large covers $\nu $ such that
$(\nu ,\nu ^{x})$ is $\varepsilon $-disjoint. We address this next.

\subsection{\label{subsec:disjointification}Disjointification lemmas}

Our methods in this section follow those in \cite{W03}. 

We first need some more notation. If $\nu ,\mu $ are covers, their
join is \[
\nu \vee \mu =\left\{ F_{n}f\, :\, F_{n}f\in \nu \textrm{ or }F_{n}f\in \mu \right\} \]
 (we use this notation instead of the more natural $\nu \cup \mu $
to avoid confusion with the set $\cup \nu $ defined above in \ref{notation:cover}).
The join of several covers $\nu _{1},\ldots ,\nu _{k}$ is denoted
by $\bigvee _{i=1}^{k}\nu _{i}$.

The \emph{restriction} of a cover $\nu $ over $E$ to a subset of
$E'\subseteq E$ is denoted by $\nu |_{E'}=\left\{ F_{n}f\in \nu \, :\, f\in E'\right\} $.
We say that $\mu $ is a \emph{subcover} of $\nu $ if $\mu \subseteq \nu $,
ie if $F_{n}f\in \nu $ whenever $F_{n}f\in \mu $. Thus if $\nu $
is a cover over $E$ and $E'\subseteq E$ we have $\nu |_{E'}\subseteq \nu $.

We also write $\min \nu =\min \left\{ n\, :\, F_{n}f\in \nu \right\} $,
and similarly $\max \nu $.

\begin{defn}
\label{notation:blowup}Let $\{F_{n}\}$ be a F\o{}lner sequence.
The \emph{blowup} of $F_{n}$ is $F_{n}^{+}=(\cup _{k<n}F_{k}^{-1})F_{n}$.
For a cover $\nu $, write $\cup ^{+}\nu =\cup _{\nu }F_{n}^{+}f$.
\end{defn}
The point of this is that if $F_{k}f$ and $F_{n}g$ are two translates
with $k<n$, then $f\notin F_{n}^{+}g$ iff $F_{k}f\cap F_{n}g=\emptyset $.
Thus if $F_{k}f$ is a set and $\nu $ some cover, and $k<\min \nu $,
then $f\notin \cup ^{+}\nu $ iff $F_{k}f\cap (\cup \nu )=\emptyset $.

For a F\o{}lner sequence $\{F_{n}\}$, the condition $|F_{n}^{+}|\leq C|F_{n}|$
is equivalent to saying that $\{F_{n}\}$ is tempered with constant
$C$.

The following properties will be useful. They state that if $\left\{ F_{n}\right\} $
is tempered then the size of $\cup ^{+}\nu $, $\cup (\nu \vee \nu ^{x})$
and $\cup ^{+}(\nu \vee \nu ^{x})$ are at most a constant times the
size of $\cup \nu $, assuming some disjointness criteria from $\nu $
or $(\nu ,\nu ^{x})$.

\begin{lem}
\label{lemma:epsilon-disjointness-means-small-blowup}Suppose $\{F_{n}\}$
is tempered with constant $C$. If $\nu $ is an $\varepsilon $-disjoint
cover then $|\cup ^{+}\nu |\leq \frac{C}{1-\varepsilon }|\cup \nu |$.
\end{lem}
\begin{proof}
We have \[
|\cup ^{+}\nu |=|\bigcup _{\nu }F_{n}^{+}f|\leq \sum _{\nu }|F_{n}^{+}f|\leq C\sum _{\nu }|F_{n}f|\leq \frac{C}{1-\varepsilon }|\cup \nu |\]
 (the last inequality follows from the fact that $\nu $ is $\varepsilon $-disjoint). 
\end{proof}
\begin{lem}
\label{lemma:temperedness-and-pair-epsilon-disjointness-means-small-blowup}If
$(\nu ,\nu ^{x})$ is $\varepsilon $-disjoint then $|\cup (\nu \vee \nu ^{x})|\leq \frac{2}{1-\varepsilon }|\cup \nu |$.
If in addition $\left\{ F_{n}\right\} $ is tempered with constant
$C$ then $|\cup ^{+}(\nu \vee \nu ^{x})|\leq \frac{2C}{1-\varepsilon }|\cup \nu |$.
\end{lem}
\begin{proof}
If $(\nu ,\nu ^{x})$ is $\varepsilon $-disjoint then it is clear
from the definitions that $\nu $ is $\varepsilon $-disjoint. We
therefore have\[
|\cup (\nu \vee \nu ^{x})|\leq \sum _{\nu }(|F_{n}f|+|(F_{n}f)^{x}|)=2\sum _{\nu }|F_{n}f|\leq \frac{2}{1-\varepsilon }|\cup \nu |.\]
 If $\left\{ F_{n}\right\} $ is tempered with constant $C$,\[
|\cup ^{+}(\nu \vee \nu ^{x})|\leq 2\sum _{\nu }|F_{n}^{+}f|\leq 2C\sum _{\nu }|F_{n}f|\leq \frac{2C}{1-\varepsilon }|\cup \nu |.\qedhere \]

\end{proof}
A cover $\nu $ over $E$ is constant if for some $n_{0}$ every member
of $\nu $ is of the form $F_{n_{0}}f$. We then say that $\nu $
is an $n_{0}$-cover. For any cover $\nu $, we can write $\nu $
as $\nu =\bigvee _{i=1}^{I}\nu _{i}$ where each $\nu _{i}$ is a
constant $i$-cover. We say that $\nu _{i}$ is the $i$-th \emph{level}
of $\nu $, and that $\{\nu _{i}\}$ is the \emph{decomposition of
$\nu $ into levels.}

Finally, we say that a cover $\nu $ over $E$ is \emph{simple} if
for each $f\in E$ there is exactly one $n$ such that $F_{n}f\in \nu $.
Any cover $\nu $ over $E$ has a simple subcover over $E$.

The first step in proving that large covers $\nu $ exist such that
$(\nu ,\nu ^{x})$ is $\varepsilon $-disjoint is the following lemma. 

\begin{lem}
\label{lemma:coding-disjointification} Let $0<\varepsilon <\frac{1}{2}$
and $x\in \Omega $. Let $\nu $ be a cover over a set $E$. Then
there exists a simple subcover $\mu \subseteq \nu $ over a set $E'=\left\{ f_{1},\ldots ,f_{I}\right\} \subseteq E$
such that $(\mu ,\mu ^{x})$ is $\varepsilon $-disjoint and $|\cup \mu |\geq \frac{\varepsilon }{8+2C}|E|$.
\end{lem}
\begin{proof}
With regards the requirement that $\mu $ be simple, we merely note
that we can replace $\nu $ with a simple subcover of $\nu $ and
proceed from there.

We first prove the lemma under the assumption that $\nu $ is constant
and then proceed to the general case.

Suppose $\nu $ is a constant $n$-cover over $E$. We will show that
there exists $E'\subseteq E$ such that for $\mu =\nu |_{E'}$ we
have $(\mu ,\mu ^{x})$ is $\varepsilon $-disjoint and $|\cup \mu |\geq \frac{\varepsilon }{8}|E|$.
Let $E'=\left\{ f_{1},\ldots ,f_{I}\right\} \subseteq E$ be a maximal
set such that for $\mu =\nu |_{E'}$, $(\mu ,\mu ^{x})$ is $\varepsilon $-disjoint.
Set $D=\cup (\mu \vee \mu ^{x})$. If $|D\cap E|>\frac{1}{2}|E|$
then by lemma \ref{lemma:temperedness-and-pair-epsilon-disjointness-means-small-blowup},
\[
|\cup \mu |\geq \frac{1-\varepsilon }{2}|D|\geq \frac{1}{4}|D|\geq \frac{1}{4}|D\cap E|>\frac{1}{8}|E|\]
 so we are done. Otherwise $|D\cap E|\leq \frac{1}{2}|E|$. Then by
minimality, for each $f\in E\setminus D$ it must hold that $|F_{n}f\cap D|\geq \varepsilon |F_{n}|$.
Since each $g\in D$ can lie in at most $|F_{n}|$ right translates
of $F_{n}$, we have $|D|\geq \frac{1}{|F_{n}|}\cdot \varepsilon |D\setminus E|\cdot |F_{n}|\geq \frac{\varepsilon }{2}|E|$,
so \[
|\cup \mu |\geq \frac{1-\varepsilon }{2}|D|\geq \frac{(1-\varepsilon )\varepsilon }{4}|E|\geq \frac{\varepsilon }{8}|E|.\]

Returning to the case where the $\nu $ is not constant, let $\nu =\bigvee _{i=1}^{I}\nu _{i}$
be the levels of $\nu $. Since $\nu $ is simple, the set of centers
of $\nu _{i}$ and $\nu _{j}$ are disjoint fro $i\neq j$. We define
inductively for $i=I$ down to $i=1$ subcovers $\mu _{i}$ of $\nu _{i}$
over sets $E'_{i}=\{f_{i1},\ldots ,f_{iM(i)}\}$ such that $(\mu _{i},\mu _{i}^{x})$
is $\varepsilon $-disjoint, as follows: Assume we have defined $\mu _{j}$
for $j>i$ and define \[
\nu '_{i}=\left\{ F_{i}f\in \nu _{i}\, :\, F_{i}f\cap \left(\cup (\mu _{j}\vee \mu _{j}^{x})\right)=\emptyset \textrm{ for every }j>i\right\} .\]
 This is equivalent to saying that $\nu '_{i}=\nu _{i}|_{E_{i}}$,
where \[
E_{i}=\dom \nu _{i}\setminus \bigcup _{j>i}\left(\cup ^{+}(\mu _{j}\vee \mu _{j}^{x})\right).\]
 Now apply the one-level case proved above to $\nu '_{i}$ to obtain
a subcover $\mu _{i}$ which is the restriction of $\nu _{i}$ to
$E'_{i}=\{f_{i1},\ldots ,f_{iM(i)}\}\subseteq E_{i}$.

Finally, set $E'=\cup E'_{i}$, and order $E'$ as in \[
E'=\{f_{i}\}_{i=1}^{M}=\{f_{I1},\ldots ,f_{IM(I)},\ldots ,f_{11},\ldots ,f_{1M(1)}\}\]
 (note that the $E'_{i}$ are pairwise disjoint), and define $\mu =\cup \mu _{i}=\nu |_{E'}$.

It is easy to verify from the construction that $(\mu ,\mu ^{x})$
is $\varepsilon $-disjoint. We must verify that $\cup \mu $ is large
enough. Write $D=\cup ^{+}(\mu \vee \mu ^{x})$ and note that $E_{i}\supseteq \dom \nu _{i}\setminus D$
so \[
\cup _{i}E_{i}\supseteq (\cup _{i}\dom \nu _{i})\setminus D=E\setminus D.\]
Hence\[
|\cup \mu |=\sum _{i=1}^{I}|\cup \mu _{i}|\]
 which by the definition of $\mu _{i}$ and the one-level case proved
at the beginning implies\begin{eqnarray*}
|\cup \mu | & \geq  & \frac{\varepsilon }{8}\sum _{i=1}^{I}|E_{i}|\\
 & \geq  & \frac{\varepsilon }{8}|\bigcup _{i=1}^{I}E_{i}|\\
 & \geq  & \frac{\varepsilon }{8}(|E|-|D|)\\
 & \geq  & \frac{\varepsilon }{8}(|E|-\frac{2C}{1-\varepsilon }|\cup \mu |)
\end{eqnarray*}
 the last inequality since by lemma \ref{lemma:temperedness-and-pair-epsilon-disjointness-means-small-blowup}
we have $|D|\leq \frac{2C}{1-\varepsilon }|\cup \mu |$. After rearranging
this is\[
|\cup \mu |\geq (\frac{8}{\varepsilon }+\frac{2C}{1-\varepsilon })^{-1}|E|\geq \frac{\varepsilon }{8+2C}|E|.\qedhere \]

\end{proof}
Next, we want to get rid of the $\varepsilon $ in the constant on
the right hand side of the last inequality. For this, the idea is
to apply the last lemma repeatedly to several {}``layers'' of covers.
If we arrange that the covers are invariant enough with respect to
each other, we can at each stage obtain a subcover with size almost
a fixed fraction of the size of what remains after the last stage,
and repeating this enough times we will get a large fraction of $E$.

\begin{defn}
\label{notation:epsilon-tempered-order-relation}For fixed $\{F_{n}\}$,
$\varepsilon >0$, if $m<n$ are integers we write $m\prec _{\varepsilon }n$
if for all $n'\geq n$ $|(\cup _{i\leq m}F_{i}^{-1})F_{n'}|\leq (1+\varepsilon )|F_{n'}|$.
For covers $\nu _{1},\nu _{2}$ we write $\nu _{1}\prec _{\varepsilon }\nu _{2}$
if $\max \nu _{1}\prec _{\varepsilon }\min \nu _{2}$.
\end{defn}
It is easy to verify that the relation $\prec _{\varepsilon }$ is
transitive. As an example of the relation $\prec _{\varepsilon }$,
note that if $\{F_{n}\}$ is $(1+\varepsilon )$-tempered, then $i\prec _{\varepsilon }(i+1)$
for all $i$. Also note that for any F\o{}lner sequence $\{F_{n}\}$
and for any $m$ and $\varepsilon $, $m\prec _{\varepsilon }n$ is
true for large enough $n$.

\begin{defn}
\label{notation:epsilon-blowup}For $n\in \mathbb{N}$ and $\varepsilon >0$,
let \[
F_{n}^{(\varepsilon )}=(\bigcup _{i:i\prec _{\varepsilon }n}F_{i}^{-1})F_{n}.\]
 For a cover $\nu $ let $\cup ^{\varepsilon }\nu =\cup \{F_{n}^{(\varepsilon )}f\, :\, F_{n}f\in \nu \}$.
\end{defn}
By definition of the $\prec _{\varepsilon }$ relation, $|F_{n}^{(\varepsilon )}|\leq (1+\varepsilon )|F_{n}|$.
The following lemmas are proved like \ref{lemma:epsilon-disjointness-means-small-blowup}
and \ref{lemma:temperedness-and-pair-epsilon-disjointness-means-small-blowup}:

\begin{lem}
\label{lemma:small-epsilon-blowup}If $\nu $ is an $\varepsilon $-disjoint
cover, then $|\cup ^{\varepsilon }\nu |\leq \frac{1+\varepsilon }{1-\varepsilon }|\cup \nu |$,
and if $(\nu ,\nu ^{x})$ is $\varepsilon $-disjoint then $|(\cup ^{\varepsilon }\nu )\cup (\cup ^{\varepsilon }\nu ^{x})|\leq \frac{2(1+\varepsilon )}{1-\varepsilon }|\cup \nu |$. 
\end{lem}
\begin{prop}
\label{lemma:iterated-coding-disjointification}Let $0<\varepsilon <\frac{1}{2}$
and $x\in \Omega $. Then there exists an integer $N$ such that for
any sequence of $N$ covers $\nu _{1}\prec _{\varepsilon }\nu _{2}\prec _{\varepsilon }\ldots \prec _{\varepsilon }\nu _{N}$
over $E$ there a exists simple subcover $\mu \subseteq \bigvee _{k=1}^{N}\nu _{k}$
such that $(\mu ,\mu ^{x})$ is $\varepsilon $-disjoint and $|\cup \mu |\geq \frac{1}{8}|E|$. 
\end{prop}
\begin{proof}
We will describe a subcover $\mu \subseteq \bigvee _{k=1}^{N}\nu _{k}$,
and show that if $N$ was large enough to begin with, the desired
properties hold.

The proof is similar to lemma \ref{lemma:coding-disjointification}.
Define inductively, from $i=N$ down to $i=1$, subcovers $\mu _{i}$
of $\nu _{i}$ defined over sets $E'_{i}=\{f_{i1},\ldots ,f_{iM(i)}\}\subseteq E$
as follows: Assuming we have defined $\mu _{j},E'_{j}$ for $j>i$,
set\[
E_{i}=E\setminus \bigcup _{j\geq k}\left(\dom \mu _{j}\, \cup \, (\cup ^{\varepsilon }(\mu _{j}\vee \mu _{j}^{x})\right)\]
 and $\nu '_{i}=\nu _{i}|_{E_{i}}$. $\mu _{i}$ and $E'_{i}$ are
obtained by and applying lemma \ref{lemma:coding-disjointification}
to $\nu '_{i}$. Finally, define $E'=\cup _{i\leq I}E'_{i}$ and order
$E'$ as in lemma \ref{lemma:coding-disjointification} \[
E'=\{f_{i}\}_{i=1}^{M}=\{f_{I1},\ldots ,f_{IM(I)},\ldots ,f_{11},\ldots ,f_{1M(1)}\}\]
 (as before the $E'_{i}$ are pairwise disjoint). Set $\mu =\cup \mu _{j}$
It is easy to verify from the construction that $\mu $ is simple
and $(\mu ,\mu ^{x})$ is $\varepsilon $-disjoint. We now show that
$|\cup \mu |\geq \frac{1}{8}|E|$.

Set $D=\cup ^{\varepsilon }(\mu \vee \mu ^{x})$. Then $|D|\leq \frac{2(1+\varepsilon )}{1-\varepsilon }|\cup \mu |$.
Also, clearly $|\dom \mu |\leq |\cup \mu |$ because $\mu $ is $\varepsilon $-disjoint.
Thus for each $i=1,\ldots ,N$,\begin{eqnarray*}
|E_{i}| & \geq  & |E|-|D|-|\dom \mu |\\
 & \geq  & |E|-(\frac{2(1+\varepsilon )}{1-\varepsilon }+1)|\cup \mu \\
 & \geq  & |E|-\frac{3-\varepsilon }{1-\varepsilon }|\cup \mu |
\end{eqnarray*}
 so by lemma \ref{lemma:coding-disjointification}\begin{eqnarray*}
|\cup \mu | & = & \sum _{i=1}^{N}|\cup \mu _{i}|\geq \\
 & \geq  & \sum _{i=1}^{N}\frac{\varepsilon }{8+2C}|E_{i}|\\
 & \geq  & \sum _{i=1}^{N}\frac{\varepsilon }{8+2C}(|E|-\frac{3-\varepsilon }{1-\varepsilon }|\cup \mu |)\\
 & \geq  & \frac{\varepsilon N}{8+2C}(|E|-\frac{3-\varepsilon }{1-\varepsilon }|\cup \mu |)
\end{eqnarray*}
 and after rearranging this is\[
|\cup \mu |\geq \left(\frac{8+2C}{\varepsilon N}+\frac{(3-\varepsilon )}{(1-\varepsilon )}\right)^{-1}|E|\]
 for $N$ large enough this is greater than $\frac{1}{8}|E|$
\end{proof}
We will also need the following lemma, which can be proved by the
same techniques as above. It was first proved by E. Lindenstrauss
in \cite{Lin01}, or see \cite{W03} for a different approach more
similar to ours.

\begin{lem}
\label{lemma:epsilon-disjointification}Let $0<\varepsilon <\frac{1}{2}$.
Then there exists an integer $N$ such that for any sequence of $N$
covers $\nu _{1}\prec _{\varepsilon }\nu _{2}\prec _{\varepsilon }\ldots \prec _{\varepsilon }\nu _{N}$
over $E$, there exists a simple subcover $\nu \subseteq \bigvee _{i=1}^{N}\nu _{i}$
such that $\nu $ is $\varepsilon $-disjoint and $|\cup \nu |\geq (1-2\varepsilon )|E|$. 
\end{lem}
We end this section with a combinatorial result akin to lemma \ref{lemma:number-of-small-subsets}:

\begin{lem}
\label{lemma:number-of-epsilon-disjoint-covers}(Number of $\varepsilon $-disjoint
covers) Let $L<N$ and $0<\varepsilon <\frac{1}{2}$. Then the number
of simple $\varepsilon $-disjoint covers $\nu $ with centers in
$F_{N}$, $\min \nu >L$ and $\cup \nu \subseteq F_{N}$ is at most
$2^{\rho (\varepsilon ,L)|F_{N}|}$, where $\lim _{L\rightarrow \infty }\rho (\varepsilon ,L)=0$,
uniformly in $\varepsilon $. 
\end{lem}
\begin{proof}
Each cover $\nu =\left\{ F_{n(i)}f_{i}\right\} _{i=1}^{I}$ of the
type we are interested in is determined by its centers $\left\{ f_{1},\ldots ,f_{I}\right\} $
and the values of $n(i)$ for $i=1,\ldots ,I$. If $\nu $ is $\varepsilon $-disjoint
and $\cup \nu \subseteq F_{N}$, we have\[
|F_{N}|\geq |\cup \nu |\geq (1-\varepsilon )\sum _{i=1}^{I}|F_{n(i)}f_{i}|\geq (1-\varepsilon )\cdot I\cdot L\]
 since $|F_{n}f|\geq n$ and $\min \nu >L$. Thus $I\leq \frac{1}{(1-\varepsilon )L}|F_{N}|$.
Since we only consider simple covers, the number of possible sets
of centers for such $\varepsilon $-disjoint covers $\nu $ with centers
in $F_{N}$ is bounded by the number of subsets of size at most $\frac{1}{(1-\varepsilon )a_{L}}|F_{N}|$
of $|F_{N}|$, which by lemma \ref{lemma:number-of-small-subsets}
is bounded above by \[
2^{\delta (\frac{1}{(1-\varepsilon )L})|F_{N}|}.\]

We now calculate how many $\varepsilon $-disjoint covers $\nu $
exist which satisfy our assumptions and have $\dom \nu =E$ for some
fixed set $E$. We use a coding argument: Given $\nu =\left\{ F_{n(i)}f_{i}\right\} _{i=1}^{I}$
as above, we write down the values $n(i)$ in binary notation in some
pre-determined order (eg fix an order on $G$ and write $n(i)$ down
according to the order of the $f_{i}$). In order to decode the resulting
concatenated sequence of $0$s and $1$s we must insert {}``punctuation''
into the resulting concatenation of numbers. We therefore code each
binary digit $0$ as $00$ and each $1$ as $11$, and terminate each
number with the pair $01$. Using this scheme, we can clearly recover
$\nu $ from $E$ and the codeword. In order to code $n(i)$ in this
manner we use $2\log n(i)+2$ bits. Recalling that we always assume
$|F_{n}|\geq n$, the length of the resulting codeword is\begin{eqnarray*}
\sum _{i=1}^{I}(2\log n(i)+2) & \leq  & \sum _{i=1}^{I}(2+2\log |F_{n(i)}f_{i}|)\\
 & \leq  & 2\cdot I+2\sum _{i=1}^{I}|F_{n(i)}|\frac{\log |F_{n(i)}|}{|F_{n(i)}|}\\
 & \leq  & \frac{2}{(1-\varepsilon )L}|F_{N}|+2\cdot \frac{\log L}{L}\sum _{i=1}^{I}|F_{n(i)}f_{i}|\\
 & \leq  & \frac{2}{(1-\varepsilon )L}|F_{N}|+\frac{2\log L}{L}\cdot \frac{|\cup F_{n(i)}f_{i}|}{1-\varepsilon }\\
 & \leq  & \frac{2}{(1-\varepsilon )L}|F_{N}|+\frac{2\log L}{(1-\varepsilon )L}|F_{N}|.
\end{eqnarray*}
 Thus there are at most $2^{(\frac{2}{(1-\varepsilon )L)}+\frac{\log L}{(1-\varepsilon )L})|F_{N}|}$
such covers over a fixed set $E$.

Putting this all together, we see that the number of $\varepsilon $-disjoint
covers $\nu $ with $\min \nu >L$ and $\cup \nu \subseteq F_{N}$
is at most\[
\exp \left(\ln 2\cdot \left(\frac{2\log L}{(1-\varepsilon )L}+\delta (\frac{1}{(1-\varepsilon )L})+\frac{2}{(1-\varepsilon )L}\right)|F_{N}|\right)\]
 which proves the lemma. 
\end{proof}

\subsection{\label{subsec:lower-bound-for-T-part-2}Proof of the lower bound
for $T_{*}$ (Part II)}

\begin{proof}
(of theorem \ref{thm:lower-bound-for-T}) Fix $\left\{ F_{n}\right\} $,$\left\{ W_{n}\right\} $
and the $\left(X_{g}\right)_{g\in G}$, defined on $(\Omega ,\mathcal{F},P)$.
Suppose $P(T_{*}^{(F,W)}<h-\varepsilon )>p>0$. We are going to use
lemma \ref{lemma:iterated-coding-disjointification} to construct
a sequence $\{c_{n}\}$ of $F_{n}$-codes such that eventually almost
surely $x(F_{n})$ can be reconstructed from $c_{n}(x)$, and \[
\frac{1}{|F_{n}|}\ell (c_{n}(x))\leq h-\frac{p}{5}\varepsilon .\]
 This would contradict the fact that entropy is a lower bound for
the achievable coding rates (section \ref{subsec:coding}), since
we can proceed as in lemma \ref{lemma:making-faithful-codes-from-nonfaithful-codes}
and turn $\left\{ c_{n}\right\} $ into a faithful code with the same
rate.

Fix a $\eta >0$ and a series $\left\{ \Phi _{n}\right\} $ of sets
of entropy-typical $F_{n}$-patterns, $\Phi _{n}\subseteq \Sigma ^{F_{n}}$,
such that for every $\varphi \in \Phi $, \[
2^{-(h+\eta )|F_{n}|}\leq P([\varphi ])\leq 2^{-(h-\eta )|F_{n}|}\]
and $x\in [\Phi _{n}]$ eventually almost surely. For $x\in \Omega $
and $N$ fixed, consider the pattern $x(F_{N})$; the codeword $c_{N}(x)$
depends only on this pattern, and is constructed as follows: 
\begin{description}
\item [Step~1]Consider the sets $F_{k}f\subseteq F_{N}$ such that $T_{k}(fx)<2^{(h-\varepsilon /2)|F_{k}|}$,
and from this collection try to extract a cover $\mu $ such that
$(\mu ,\mu ^{x})$ is $\eta $-disjoint and very large. To be precise,
we require that $|\cup \mu |\geq \frac{p}{8}|F_{N}|$. However, we
try to choose $\mu $ so that the next step is possible:
\item [Step~2]Consider the sets $F_{k}f\subseteq F_{N}$ for which $fx\in [\Phi _{k}]$,
and from this collection try to extract an $\eta $-disjoint cover
$\nu $ which is also disjoint from $\cup \mu $. We choose $\nu $
as large as possible; we require that $\cup \nu $ and $\cup \mu $
together cover all but a $3\eta $-fraction of $F_{N}$.
\item [Step~3]If we cannot find $\mu ,\nu $ as in steps 1,2 with the
required sizes, we define $c_{N}(x)$ to be the empty word. 
\item [Step~4]Otherwise, let $A_{1}$ be the set guarenteed by lemma \ref{lemma:pair-decoding}
when applied to $\mu $, so $A_{1}\subseteq \cup \mu $ and $|A_{1}|<2\eta |\cup \mu |\leq 2\eta |F_{N}|$,
and $A_{1}$ has the property that if pattern $x_{1}(A\cup (\cup \mu ^{x}))$
is known and if the elements $f^{x,k}$ are known for every $F_{k}f\in \mu $
then $x(\cup \mu )$ can be deduced. Also, let $A_{2}=F_{N}\setminus ((\cup \mu )\cup (\cup \nu ))$.
\item [Step~5]Now define (assuming steps 1 and 2 were successful) $c_{n}(x)=\gamma _{1}\gamma _{2}\gamma _{3}$,
where

\begin{itemize}
\item $\gamma _{1}$ encodes the pattern $x(A_{1}\cup A_{2})$.
\item $\gamma _{2}$ encodes the pattern $x(\cup \nu )$.
\item $\gamma _{3}$ encodes the cover $\mu $ along with the elements $f^{x,k}$
for each $F_{k}f\in \mu $.
\end{itemize}
\end{description}
It is clear from the definition of $A_{1}$ that if steps 1,2 succeeded
then we can reconstruct $x(F_{N})$ from $c_{N}(x)$. It remains to
show that for almost every $x$, steps 1 and 2 will succeed for large
enough $N$, and that in this case the encoding, as described in step
5, can be carried out in such a way that $\ell (c_{N}(x))\leq (h-p/10)|F_{N}|$. 

We begin by showing that for large enough $N$, steps 1 and 2 will
succeed. Fix integers $M,L$ to be determined later, and select a
sequence of intervals $\{[a_{i};b_{i}]\}_{i=1}^{2L}$ such that $a_{1}=M$,
$b_{i}\prec _{\eta }a_{i+1}$ and 
\begin{enumerate}
\item For a set of $x\in \Omega $ with probability greater than $p$ there
exist indices $m_{i}(x)\in [a_{i};b_{i}]$ such that $T_{m_{i}(x)}(x)\leq h-\varepsilon $
for each $i$.
\item For a set of $x\in \Omega $ with probability greater than $1-\eta $
there exist indices $n_{i}(x)\in [a_{i};b_{i}]$ such that $x\in \Phi _{n_{i}(x)}$
for each $i$.
\end{enumerate}
Such a sequence is easy to construct using the ergodic theorem, the
fact that $P(T_{*}\leq h-\varepsilon )>p$, and our choice of $\left\{ \Phi _{n}\right\} $;
we omit the details.

From the pointwise ergodic theorem and the asymptotic invariance of
$\{F_{n}\}$, for almost every $x\in \Omega $ it holds for large
enough $N$ that\begin{equation}
\frac{1}{|F_{N}|}\#\left\{ f\in F_{N}\, \left|\begin{array}{c}
 \textrm{for }L+1\leq i\leq 2L\textrm{ there exist indices}\\
 m_{i}(f)\in [a_{i};b_{i}]\textrm{ with }T_{m_{i}(f)}(fx)\leq h-\varepsilon \textrm{ }\\
 F_{m_{i}(f)}f\subseteq F_{N}\textrm{ and }(F_{m_{i}(f)}f)^{x}\subseteq F_{N}\end{array}
\right.\right\} >p\label{eq:large-fraction-of-nearby-repetitions}\end{equation}
\begin{equation}
\frac{1}{|F_{N}|}\#\left\{ f\in F_{N}\, \left|\begin{array}{c}
 \textrm{for }1\leq i\leq L\textrm{ there exist indices}\\
 n_{i}(f)\in [a_{i};b_{i}]\textrm{ with}\\
 fx\in \Phi _{n_{i}(f)}\textrm{ and }F_{n_{i}(f)}f\subseteq F_{N}\end{array}
\right.\right\} >1-\eta .\label{eq:large-fraction-of-entropy-typical-sets}\end{equation}
 For $x,N$ such that the above holds, let $m_{i}(f)$ and $n_{i}(f)$
be as in (\ref{eq:large-fraction-of-nearby-repetitions}) and (\ref{eq:large-fraction-of-entropy-typical-sets}).

Let $E_{1}\subseteq F_{N}$ be the set described in \ref{eq:large-fraction-of-nearby-repetitions},
and $\mu _{i}=\left\{ F_{m_{i}(f)}f\right\} _{f\in E_{1}}$, so we
have covers $\mu _{i}$ over $E_{1}$. Note that $\mu _{1}\prec _{\eta }\mu _{2}\prec _{\eta }\ldots \prec _{\eta }\mu _{L}$
by our choice of $[a_{i};b_{i}]$; so assuming $L$ was chosen large
enough we can apply lemma \ref{lemma:epsilon-disjointification} to
$\mu _{1},\ldots ,\mu _{L}$ and obtain a simple subcover $\mu \subseteq \bigvee _{i=1}^{L}\mu _{i}$
over some set $E'_{1}=\{f_{1},\ldots ,f_{I}\}\subseteq E_{1}$ such
that $(\mu ,\mu ^{x})$ is $\eta $-disjoint and $|\cup \mu |\geq \frac{1}{8}|E_{1}|$,
which means that \[
|\cup \mu |\geq \frac{p}{8}|F_{N}|.\]

Let $E_{2}\subseteq F_{N}$ be the set described in \ref{eq:large-fraction-of-entropy-typical-sets},
and $\nu _{i}=\left\{ F_{n_{i}(f)}f\right\} _{f\in E_{2}}$, so we
have covers $\nu _{i}$ over $E_{2}$. As with the $\mu _{i}$, we
have $\nu _{1}\prec _{\eta }\nu _{2}\prec _{\eta }\ldots \prec _{\eta }\nu _{L}$,
and furthermore $\nu _{L}\prec _{\eta }\mu $. Let $\nu '_{i}$ be
the restriction of $\nu _{i}$ to $E_{2}\setminus \cup ^{\eta }\mu $,
so each member set of $\nu '_{i}$ is disjoint from $\cup \mu $.
Assuming $L$ is large enough we apply lemma \ref{lemma:epsilon-disjointification}
and obtain an $\eta $-disjoint simple subcover $\nu \subseteq \bigvee _{i=1}^{L}\nu '_{i}$
over some set $E'_{2}\subseteq E_{2}\setminus (\cup ^{\eta }\mu )$
such that $|\cup \nu |\geq (1-2\eta )|E_{2}\setminus (\cup ^{\eta }\mu )|$.
We then have \[
|\cup \nu |\geq (1-2\eta )\left(|E_{2}|-|\cup ^{\eta }\mu |\right)\geq (1-2\eta )\left((1-\eta )|F_{N}|-\frac{1+\eta }{1-\eta }|\cup \mu |\right).\]

By construction $\mu ,\nu $ are disjoint, and putting the last two
paragraphs together, we have \[
|\cup (\mu \vee \nu )|=|\cup \mu |+|\cup \nu |\geq \]
 \[
\geq (1-(1-2\eta )\frac{1+\eta }{1-\eta })|\cup \mu |+(1-2\eta )(1-\eta )|F_{N}|\geq (1-3\eta )|F_{N}|.\]
 Thus all but a $3\eta $-fraction of $F_{N}$ is covered by $\mu \cup \nu $,
and we have shown that for almost every $x$, for large enough $N$
there exist $\mu ,\nu $ as required by steps 1 and 2 of the coding
construction. 

We complete the proof by giving a more detailed description of the
encoding of $c_{N}(x)$ (in the case where steps 1,2 succeeded) by
specifying exactly how to construct $\gamma _{1},\gamma _{2},\gamma _{3}$,
and show that the codeword length is indeed bounded by $(h-p/5)|F_{N}|$.

Recall the functions $\delta (\cdot )$ and $\rho (\cdot ,\cdot )$
defined in lemmas \ref{lemma:number-of-small-subsets} and \ref{lemma:number-of-epsilon-disjoint-covers}
respectively. In the estimates below we ignore rounding errors; the
reader may verify that these may be taken into account without disrupting
the proof.

Fix $x,N$ so that steps 1,2 of the construction succeeded. Let $\mu ,\nu ,A_{1},A_{2}$
be as in the construction.

Write $H=A_{1}\cup A_{2}$. The word $\gamma _{1}$ will first encode
the set $H$ and then the values $x(f)$ for $f\in H$. The former
takes $\delta (5\eta )|F_{N}|$ bits, because by lemma \ref{lemma:number-of-small-subsets}
we can costruct a list of all $2^{\delta (5\eta )|F_{n}|}$ subsets
of $F_{n}$ with size at most $5\eta |F_{n}|$ and then describe $H$
by giving its index in the list. Next, using some fixed ordering of
$G$ (which induces an order on $H$), we write down the symbols $x(f)$
for $f\in H$. For this we need an additional $|H|\log \Sigma \leq (5\eta \log \Sigma )|F_{N}|$
bits. So all together $x(H)$ can be coded in $(\delta (5\eta )+5\eta \log \Sigma )|F_{N}|$
bits.

The word $\gamma _{2}$ will be coded by specifying first the cover
$\nu $ and then the values $x(F_{k}f)$ for $F_{k}f\in \nu $. The
first task is accomplished using $\rho (\eta ,M)|F_{N}|$ bits by
lemma \ref{lemma:number-of-epsilon-disjoint-covers}. Now for each
$F_{k}f\in \nu $, according to some fixed order, we record the pattern
$x(F_{k}f)$ by giving the index of the pattern $x(F_{k}f)$ in $\Phi _{k}$.
This requires $(h+\eta )|F_{k}|$ bits by our choice of $\{\Phi _{n}\}$,
because $|\Phi _{k}|\leq 2^{(h+\eta )|F_{k}|}$, and because by the
definition of $\nu $, $fx\in [\Phi _{k}]$ whenever $F_{k}f\in \nu $.
Thus to describe all the patterns $x(F_{k}f)$ for $F_{k}f\in \nu $,
we need\[
\sum _{\nu }(h+\eta )|F_{k}f|=(h+\eta )\sum _{\nu }|F_{k}f|\leq \frac{h+\eta }{1-\eta }|\cup \nu |\]
 bits. All together, we see that we can code $\gamma _{2}$ with at
most $\rho (\eta ,M)|F_{N}|+\frac{h+\eta }{1-\eta }|\cup \nu |$ bits.

The word $\gamma _{3}$ is coded by first describing $\mu $ and then,
for each $F_{k}f\in \mu $, we describe $f^{x,k}$. The former is
again achieved with $\rho (\eta ,M)|F_{N}|$ bits (the same calculation
as above applies), while the latter requires at most $(h-\varepsilon )|F_{\mu (f)}|$
bits for each $F_{k}f\in \mu $, by definition of $f^{x,k}$. Note
that this is where we use the assumption that $\{W_{n}\}$ is increasing:
for if it were not, then in order to record $f^{x,k}$ we would need
to specify also which set $W_{n}$ it is from. In other words, we
would need to record $R_{k}(fx)$, and this could be expensive (there
is no upper bound on $R_{k}(fx)$). However, since $\{W_{n}\}$ is
increasing, we can assume that $f^{x,k}$ comes from the largest set
$W_{n}$ satisfying $|W_{n}|<2^{(h-\varepsilon )|F_{k}|}$). Thus
estimating as we did for $\gamma _{2}$ we find that $\gamma _{3}$
can be coded in \[
\rho (\eta ,M)|F_{N}|+(h-\varepsilon )\sum _{\mu }|F_{k}f|\leq \rho (\eta ,M)|F_{N}|+\frac{h-\varepsilon }{1-\eta }|\cup \mu |\]
 bits.

Putting all this together, we find\begin{eqnarray*}
\ell (c_{N}(x)) & = & \ell (\gamma _{1})+\ell (\gamma _{2})+\ell (\gamma _{3})\\
 & \leq  & \left\{ \begin{array}{c}
 \left(\delta (7\eta )+5\eta \log \Sigma \right)|F_{N}|+\\
 +\left(\rho (\eta ,M)|F_{N}|+\frac{h+\eta }{1-\eta }|\cup \nu |\right)\\
 +\left(\rho (\eta ,M)|F_{N}|+\frac{h-\varepsilon }{1-\eta }|\cup \mu |\right)\end{array}
\right\} \\
 & \leq  & \left\{ \begin{array}{c}
 \left(\delta (7\eta )+5\eta \log \Sigma +2\rho (\eta ,M)\right)|F_{N}|+\\
 +\frac{h+\eta }{1-\eta }\left|(\cup \nu )\cup (\cup \mu )\right|-\frac{\varepsilon -\eta }{1-\eta }|\cup \mu |\end{array}
\right\} \\
 & \leq  & (h+\tau (\eta ,M))|F_{N}|-\frac{\varepsilon -\eta }{1-\eta }\cdot \frac{p}{8}|F_{N}|
\end{eqnarray*}
 where $\tau (\eta ,M)\rightarrow 0$ as $M\rightarrow \infty $ uniformly
in $\eta $. So for $M$ large enough and $\eta $ small enough we
have\[
\ell (c_{N}(x))\leq (h-\frac{p}{10}\varepsilon )|F_{N}|.\]
This completes the proof.
\end{proof}

\subsection{\label{subsec:lower-bound-for-U}The lower bound for \protect$U_{*}$}

In this section we prove that

\[
U_{*}^{(F,W)}(x)=\liminf _{k\rightarrow \infty }\, \inf _{n}-\frac{1}{|F_{k}|}\log (U_{k,n}^{(F,W)}(x)-\frac{1}{|W_{n}|})\geq h\]
 almost surely if $\left\{ F_{n}\right\} $,$\left\{ W_{n}\right\} $
satisfy the conditions of theorem \ref{thm:intro-lower-bound-for-U}.

As was already mentioned in the introduction, the correction term
$-\frac{1}{|W_{n}|}$ in the definition of $U_{*}$ is necessary because
$U_{k,n}$ counts the central pattern $x(F_{k})$; without it, if
for instance $1_{G}\in W_{1}$, we would have for every $k$ that
$U_{k,1}\geq \frac{1}{|W_{1}|}$, and this implies $U_{*}=0$ for
any process.

In the previous section, we defined the quantities $V_{k,n}$ and
$V^{*}$, and showed that $U^{*}=V^{*}$. We may analogously define
\[
V_{*}=\liminf _{k}\inf _{n}-\frac{1}{|F_{k}|}\log (V_{k,n}-\frac{1}{|W_{n}|})\]
 as a version of $U_{*}$ which counts all repetitions of $x(F_{k})$
in $x(W_{n})$, but we cannot, as was the case with $V^{*}$, show
that $U_{*}=V_{*}$, because when the correction factor is taken into
account, if we compare $U_{k,n}-\frac{1}{|W_{n}|}$ with $V_{k,n}-\frac{1}{|W_{n}|}$
as in equation \ref{eq:U-and-U-bar-equivalence}, we find that all
we can say is \[
U_{k,n}-\frac{1}{|W_{n}|}\leq V_{k,n}-\frac{1}{|W_{n}|}\leq |F_{k}|^{2}(U_{k,n}-\frac{1}{|W_{n}|})+\frac{|F_{k}|^{2}}{|W_{n}|}\]
 and this gives us $U_{*}\geq V_{*}$ but not the reverse inequality.

Furthermore, the fact that $U_{k,n}$ does not count repetitions of
$x(F_{k})$ which intersect $F_{k}$ is necessary for our proof that
$U_{*}\geq h$; the proof cannot be adapted to the case of $V_{*}$.
The importance of the requirement that $U_{k,n}$ not count the repetition
which intersect the original pattern stems from the following observation,
which is used in an essential way in the proof: For almost every $x$,
for large enough $k$, $U_{k,n}(x)=0$ unless $|W_{n}|>2^{(h-\varepsilon )|F_{k}|}$.
This follows from the fact that if $U_{k,n}(x)>1$ then there is within
$W_{n}$ a repetition of $x(F_{k})$ \emph{disjoint} from $F_{k}$,
and so $R_{k}(x)\leq n$, which in turn implies $T_{k}(x)\leq |W_{n}|$
(we will assume that $\{W_{i}\}$ is increasing); now use the fact
that for large enough $k$, by theorem \ref{thm:lower-bound-for-T},
$T_{k}(x)>2^{(h-\varepsilon )|F_{k}|}$. In the case of $V_{k,n}$,
if $V_{k,n}>1$ it may be because of a repetition of $x(F_{k})$ which
intersects $F_{k}$, and this argument fails. Thus we can't allow
repetitions which intersect $F_{k}$ to be counted in $U_{k,n}$ because
$R_{k}$ was defined as the first repetition of $x(F_{k})$ \emph{disjoint}
from $F_{k}$. If, however, we were able to prove an analogue of $T_{*}\geq h$
(theorem \ref{thm:lower-bound-for-T}) for the case where $R_{k}$
counts \emph{all} repetitions in $W_{n}$ of $x(F_{k})$, including
those intersecting $F_{k}$, then the proof below would work just
as well for $V_{*}$. This version of theorem \ref{thm:lower-bound-for-T}
is true, for instance, in the case of $W_{n}=[-k_{n};k_{n}]^{d}\subseteq \mathbb{Z}^{d}$.
This was proved in \cite{OW02}. 

One can, however, define\[
\hat{V}_{k,n}(x)=\frac{1}{|W_{n}|}\max \left\{ |E|\, \left|\begin{array}{c}
 1_{G}\in E\subseteq W_{n}\textrm{ and if }1_{G}\neq f\in E\textrm{ then }\\
 F_{k}\cap F_{k}f=\emptyset \, ,\, \textrm{and }fx\in [x(F_{k})]\end{array}
\right.\right\} .\]
 Then we have that \[
U_{k,n}-\frac{1}{|W_{n}|}\leq \hat{V}_{k,n}-\frac{1}{|W_{n}|}\leq |F_{k}|^{2}(U_{k,n}-\frac{1}{|W_{n}|})\]
 and so $\hat{V}_{*}=\liminf _{k\rightarrow \infty }\, \inf _{n}-\frac{1}{|F_{k}|}\log (\hat{V}_{k,n}(x)-\frac{1}{|W_{n}|})=U_{*}$.

Our proof that $U_{*}\geq h$ is based on the bound $T_{*}\geq h$.
It is not hard to see that $U_{*}^{(F,W)}\geq h$ implies $T_{*}^{(F,W)}\geq h$.
We would like to show the reverse, namely, that $U_{*}^{(F,W)}<h$
implies $T_{*}^{(F,W)}<h$. The idea behind the proof is that from
the relation $U_{k,n}(x)-\frac{1}{|W_{n}|}>2^{-(h-\varepsilon )|F_{k}|}$,
which means that there are {}``too many'' repetitions of $x(F_{k})$
in $W_{n}$, we can sometimes deduce that many of the points $g\in W_{n}$
at which $x(F_{k})$ repeats are such that another repetition of the
same pattern appears {}``too close'' to $g$. If we could show this,
we can appeal to the fact that $T_{*}\geq h$ to obtain a contradiction.
One problem here is that the window sets $W_{n}$ may be unsuitable
for capturing the idea of {}``too close''. This will happen if the
sequence $\{W_{n}\}$ grows too quickly, and in this case we will
not be able to obtain a contradiction through $T_{k}^{(F,W)}$. We
therefore introduce a second window sequence $\{Y_{n}\}$, which will
grow slowly enough that with respect to $\{Y_{n}\}$, a drop of $T_{k}^{(F,Y)}$
below entropy can be observed. We will show that for a nice enough
sequence $\{Y_{n}\}$, if with positive probability $U_{*}^{(F,W)}<h-\varepsilon $
is true, then it is not almost surely true that $T_{*}^{(F,Y)}\geq h$,
which is impossible.

We recall the following definitions, which were given in the introduction:

\begin{defn*}
Let $\{F_{n}\},\{W_{n}\}$ be sequences of finite subsets of $G$.
An increasing sequence $\{Y_{n}\}$ of finite subsets of $G$ is called
an \emph{interpolation sequence} \emph{for} $\{F_{n}\},\{W_{n}\}$
if 
\begin{enumerate}
\item $\left\{ Y_{n}\right\} $ is filling and increasing .
\item If $|W_{n}|\geq |Y_{m}|$ then $|Y_{m}W_{n}|\leq C|W_{n}|$ for some
constant $C$. 
\item For every pair of real numbers $0\leq \alpha <\beta $, for every
large enough $n$ there is an index $k$ such that $2^{\alpha |F_{n}|}\leq |Y_{k}|\leq 2^{\beta |F_{n}|}$.
\end{enumerate}
\end{defn*}
The main lemma we will need about interpolation sequences is that
too frequent recurrence of $x(F_{k})$ in $x(W_{n})$ implies that
$T_{k}^{(F,Y)}(fx)$ is too small for many $f\in W_{n}$.

\begin{lem}
\label{lemma:interpolation}Let $Y$ be $1/C$-filling and $E\subseteq W\subseteq G$
finite sets. Suppose that $|YW|\leq C|W|$. If $\frac{|E|}{|W|}\geq \alpha $
then \[
\frac{1}{|W|}|\left\{ f\in E\, :\, |Yf\cap E|>\beta \right\} |\geq \alpha -C^{2}\frac{\beta }{|Y|}.\]

\end{lem}
\begin{proof}
Let $E_{-}=\left\{ f\in E\, :\, |Yf\cap E|\leq \beta \right\} $,
and let $\left\{ Yf_{i}\right\} _{i=1}^{I}$ be an incremental cover
of $E_{-}$ with $f_{i}\in E_{-}$. Then by the basic lemma \ref{lemma:basic-filling-lemma}
we have that\[
\frac{|E_{-}|}{|\cup _{i}Yf_{i}|}\leq C\frac{\beta }{|Y|}\]
 since $\cup _{i}Yf_{i}\subseteq YW$ we have \[
|E_{-}|\leq C\frac{\beta }{|Y|}|YW|\leq C^{2}\frac{\beta }{|Y|}|W|\]
 so\[
\frac{1}{|W|}|\left\{ f\in E\, :\, \frac{|Yf\cap E|}{|Y|}>\beta \right\} |=\frac{1}{|W|}(|E|-|E_{-}|)\geq \alpha -C^{2}\frac{\beta }{|Y|}.\qedhere \]

\end{proof}
\begin{cor}
\label{lemma:high-density-,eans-too-soon-return}Let $0<\delta <\varepsilon /2$,
and let $\left\{ Y_{n}\right\} $ be an interpolation sequence for
$\left\{ F_{n}\right\} ,\left\{ W_{n}\right\} $. If for some $k,n$
and $x\in \Omega $ we have $U_{k,n}^{(F,W)}(x)\geq 2^{(h-\varepsilon )|F_{k}|}$,
and if $|W_{n}|>2^{(h-\delta )|F_{k}|}$, then if $k$ is large enough,\[
\frac{1}{|W_{n}|}|\left\{ f\in W_{n}\, \left|\begin{array}{c}
 fx\in [x(F_{k})]\, \textrm{and}\, \\
 T_{k}^{(F,Y)}(fx)<2^{(h-\varepsilon /2)|F_{k}|}\end{array}
\right.\right\} |>\frac{1}{2}2^{-(h-\varepsilon )|F_{k}|}.\]

\end{cor}
\begin{proof}
Set $E=\left\{ f\in W_{n}\, :\, fx\in [x(F_{k})]\right\} $ and $\beta =|F_{k}|^{2}$.
For $k$ large enough there exists an $m$ such that \[
2^{(h-2\varepsilon /3)|F_{k}|}\leq |Y_{m}|\leq 2^{(h-\varepsilon /2)|F_{k}|}\]
 Since $|W_{n}|>2^{(h-\delta )|F_{k}|}$, we have $|Y_{m}W_{n}|\leq C|W_{n}|$,
so by the lemma (with $\alpha =2^{-(h-\varepsilon )|F_{k}|}$)\[
\frac{1}{|W_{n}|}|\left\{ f\in E\, :\, |Y_{m}f\cap E|>|F_{k}^{2}|\right\} |\geq \; \; \; \; \; \; \; \; \; \; \; \; \; \; \; \; \; \; \; \; \; \; \; \; \; \; \; \; \; \; \]
\[
\geq 2^{-(h-\varepsilon )|F_{k}|}-C^{2}|F_{k}|^{2}2^{-(h-2\varepsilon /3)|F_{k}|}\geq \frac{1}{2}2^{-(h-\varepsilon )|F_{k}|}\]
 for $k$ large enough. Now note that if for some $f\in E$ we have
$|Y_{m}f\cap E|>|F_{k}|^{2}$ then at least one of the repetitions
of $fx(F_{k})$ within $fx(Y_{m})$ is disjoint from $fx(F_{k})$,
so $T_{k}^{(F,Y)}(fx)\leq |Y_{m}|<2^{(h-\delta )|F_{k}|}$. This completes
the proof. 
\end{proof}
The second requirement we make of $\left\{ W_{n}\right\} $ is that
it be incompressible. We recall the definition:

\begin{defn*}
An increasing sequence $\left\{ W_{n}\right\} $ of finite subsets
of $G$ with $1_{G}\in W_{n}$ is said to be \emph{incompressible
with constant $C$} (or $C$-incompressible) if for any incremental
sequence $\left\{ W_{n(i)}f_{i}\right\} $, the number of the sets
$W_{n(i)}f_{i}$ containing $1_{G}$ is at most $C$. 
\end{defn*}
An incompressible sequence $\{W_{n}\}$ is filling; if $\{W_{n}\}$
is incompressible with constant $C$ then it will be filling with
constant $\frac{1}{C}$, since for a collection $\{W_{n(i)}f_{i}\}_{i\in I}$
as above, each $g\in \cup _{i\in I}F_{n(i)}f_{i}$ belongs to at most
$C$ of the sets $W_{n(i)}f_{i}$. Incompressible sequences are superior
to filling sequences because of the following observation: If $\{W_{n}\}$
is incompressible and $\{W_{n(i)}f_{i}\}_{i=1}^{I}$ is incremental,
then for any finite $A\subseteq G$, we have \[
|A\cap \bigcup _{i\in I}W_{n(i)}f_{i}|\geq \frac{1}{C}\sum _{1\leq i\leq I}|A\cap W_{n(i)}f_{i}|.\]
 This is because every $a\in A$ is counted at most $C$ times in
the sum on the right.

The property of incompressibility is equivalent to the following property,
which may be described as being filling relative to arbitrary subsets
of $G$:

\begin{lem}
\label{lemma:incompressibility}Let $\left\{ W_{n}\right\} $ be incompressible,
$H_{1},H_{2}\subseteq G$ and suppose $\left\{ W_{n(i)}f_{i}\right\} _{i=1}^{I}$
is an incremental sequence such that for each $i=1,\ldots ,I$ we
have that \[
\alpha \leq \frac{|H_{1}\cap W_{n(i)}f_{i}|}{|H_{2}\cap W_{n(i)}f_{i}|}\leq \beta .\]
 Then \[
\frac{\alpha }{C}\leq \frac{|H_{1}\cap (\cup _{i}W_{n(i)}f_{i})|}{|H_{2}\cap (\cup _{i}W_{n(i)}f_{i})|}<C\beta \]
 (this remains true for $\alpha =0$ or $\beta =\infty $) .
\end{lem}
It is not difficult to see that if $\left\{ W_{n}\right\} $ satisfies
the conclusion of the lemma then it is incompressible.

\begin{proof}
We show for instance the lower bound.\[
\frac{|H_{1}\cap (\cup _{i}W_{n(i)}f_{i})|}{|H_{2}\cap (\cup _{i}W_{n(i)}f_{i})|}\geq \frac{\frac{1}{C}\sum |H_{1}\cap W_{n(i)}f_{i}|}{\sum |H_{2}\cap W_{n(i)}f_{i}|}\geq \]
 (we use here the fact that $\left\{ W_{n}\right\} $ is $\frac{1}{C}$-filling)\[
\geq \frac{1}{C}\sum _{i}\frac{|H_{2}\cap W_{n(i)}f_{i}|}{\sum _{j}|H_{2}\cap W_{n(i)}f_{j}|}\cdot \frac{|H_{1}\cap W_{n(i)}f_{i}|}{|H_{2}\cap W_{n(i)}f_{i}|}\geq \frac{\alpha }{C}.\qedhere \]

\end{proof}
\begin{thm}
Let $F_{n}$ be a tempered F\o{}lner sequence, $\left\{ X_{g}\right\} _{g\in G}$
an ergodic process with entropy $h$. Let $\left\{ W_{n}\right\} $
be an incompressible sequence. If there exists an interpolation sequence
for $\left\{ F_{n}\right\} ,\left\{ W_{n}\right\} $. Then \[
\liminf _{k\rightarrow \infty }\, \inf _{n}-\frac{1}{|F_{k}|}\log U_{k,n}^{(F,W)}(x)\geq h.\]

\end{thm}
\begin{proof}
Suppose the theorem is false. Then there is a set $B\subseteq \Omega $
with $P(B)>p>0$ such that for every $x\in B$,\[
\limsup _{k\rightarrow \infty }\, \sup _{n}-\frac{1}{|F_{k}|}\log U_{k,n}^{(F,W)}(x)\leq h-\varepsilon \]
 for some $\varepsilon >0$. Let $\left\{ \Phi _{n}\right\} $ be
a sequence of sets of words, $\Phi _{n}\subseteq \Sigma ^{F_{n}}$,
such that \[
\frac{1}{|F_{n}|}\log |\Phi _{n}|\rightarrow h\]
 and $x\in [\Phi _{n}]$ eventually almost surely. Fix $L$ very large;
we may assume that $x\in [\Phi _{n}]$ for every for $n\geq L$ and
$x\in B$. We may also assume that $L$ was chosen large enough that
$|\Phi _{k}|\leq 2^{(h+\delta )|F_{k}|}$ for $k>L$.

For every $x\in B$ there exists an index $k(x)>L$ and an index $n(x)$
such that \[
U_{k(x),n(x)}(x)\geq 2^{-(h-\varepsilon )|F_{k(x)}|}\]
 and without loss of generality we may assume that for every $x\in B$
we have $n(x)\leq M$ for some integer $M$. We further assume, increasing
$L$ if necessary, that for every $x\in B$, $|W_{n(x)}|>2^{(h-\varepsilon /4)|F_{k(x)}|}$,
because $T_{*}^{(F,W)}\geq h$ almost surely.

Let $x\in \Omega $ be typical in the sense that the ergodic theorem
holds for $B$ and every $[\Phi _{k}]$, $L\leq k\leq M$. By the
ergodic theorem, we can choose $N$ large enough so that\[
\frac{1}{|F_{N}|}|\left\{ f\in F_{N}\, :\, fx\in B\right\} |>p.\]
 For $f\in F_{N}$ such that $fx\in B$, we will write for brevity
$k(f)=k(fx)$ and $n(f)=n(fx)$.

\noindent \textbf{Step 1:} Fix $L\leq k\leq M$ and $\varphi \in \Phi _{k}$,
and let \[
H=\left\{ f\in F_{N}\, :\, fx\in [\varphi ]\right\} \; ,\; E=\left\{ f\in H\, :\, fx\in B\, ,\, k(f)=k\right\} .\]
 If $f\in E$ then $U_{k,n(f)}^{(F,W)}(fx)\geq 2^{-(h-\varepsilon )|F_{k}|}$.
Some of the points $g\in W_{n(f)}f$ at which $\varphi $ repeats
will belong to $E$ as well, and some won't. Write $E^{+}$ for the
set of $f\in E$ for which {}``many'' of the repetitions of $\varphi $
in $W_{n(f)}f$ are in $E$: \[
E^{+}=\left\{ f\in E\, :\, \frac{1}{|W_{n}|}|W_{n(f)}f\cap E|\geq 2^{-(h-\varepsilon /2)|F_{k}|}\right\} \]
 and\[
E^{-}=E\setminus E^{+}=\left\{ f\in E\, :\, \frac{1}{|W_{n}|}|W_{n(f)}f\cap E|<2^{-(h-\varepsilon /2)|F_{k}|}\right\} .\]
 We are interested in the elements of $F_{N}$ for which $T_{k}^{(F,Y)}<2^{-(h+\varepsilon /4)|F_{k}|}$.
Let\[
D=\left\{ f\in H\, :\, T_{k}^{(F,Y)}(fx)<2^{-(h+\varepsilon /2)|F_{k}|}\right\} .\]

Now if $f\in E^{+}$, then by lemma \ref{lemma:high-density-,eans-too-soon-return}
we have that\[
|(E\cap W_{n(f)}f)\cap D|\geq \frac{1}{2}|E\cap W_{n(f)}f|.\]
 Select an incremental cover $\left\{ W_{n(f_{i})}f_{i}\right\} _{i=1}^{I}$
of $E^{+}$ with $f_{i}\in E^{+}$. For each $i$ we have that the
relative density of $D\cap E$ in $W_{n(f_{i})}f_{i}$ is at least
half the density of $E$, so from lemma \ref{lemma:incompressibility}
we see that \begin{equation}
\frac{|(E\cap D)\cap (\cup W_{n(f_{i})}f_{i})|}{|E\cap (\cup W_{n(f_{i})}f_{i})|}\geq \frac{1}{2C}\label{eq:First-Use-Of-Incompressibility}\end{equation}
 and therefore, since $E^{+}\subseteq \cup W_{n(f_{i})}f_{i}$, we
have\begin{equation}
|E\cap D|\geq \frac{1}{2C}|E^{+}|.\label{eq:U-lower-bound-step-1-a}\end{equation}

On the other hand, if $f\in E^{-}$, we have that the number of elements
of $H$ in $W_{n(f)}f$ is at least $2^{(\varepsilon /2)|F_{k}|}$
times the number of elements of $E$ in $W_{n(f)}f$, so \[
|W_{n(f)}f\cap D|\geq \frac{1}{2}|W_{n(f)}f\cap H|\geq \frac{1}{2}2^{(\varepsilon /2)|F_{k}|}\cdot |E\cap W_{n(f)}f|.\]
 Now select an incremental cover $\left\{ W_{n(f_{i})}f_{i}\right\} _{i=1}^{I}$
of $E^{-}$. From lemma we have\begin{equation}
\frac{|D\cap (\cup W_{n(f_{i})}f_{i})|}{|E\cap (\cup W_{n(f_{i})}f_{i})|}\geq \frac{1}{2C}2^{(\varepsilon /2)|F_{k}|}\label{eq:Second-Use-Of-Incompressibility}\end{equation}
 and therefore, since $E^{-}\subseteq \cup W_{n(f_{i})}f_{i}$, we
have \begin{equation}
|D|\geq \frac{1}{2C}2^{(\varepsilon /2)|F_{k}|}|E^{-}|.\label{eq:U-lower-bound-step-1-b}\end{equation}

We now face two alternatives: 
\begin{enumerate}
\item If $|E^{+}|\geq \frac{1}{2}|E|$ then equation \ref{eq:U-lower-bound-step-1-a}
gives us\[
|E\cap D|\geq \frac{1}{2C}|E^{+}|\geq \frac{1}{4C}|E|.\]

\item Otherwise, if $|E_{k}^{\varphi ,-}|>\frac{1}{2}|E_{k}^{\varphi }|$,
equation \ref{eq:U-lower-bound-step-1-b} gives \[
|D|\geq \frac{1}{2C}2^{(\varepsilon /2)|F_{k}|}|E^{-}|\geq \frac{1}{4C}2^{(\varepsilon /2)|F_{k}|}|E|.\]

\end{enumerate}
\noindent \textbf{Step 2:} For fixed $k$ and $\varphi \in \Sigma ^{F_{k}}$
write $D^{\varphi },E^{\varphi }$ for the sets $D,E$ of step 1,
thus making explicit the dependence on $\varphi $ which was previously
suppressed. Note that the collections $\left\{ D^{\varphi }\right\} _{\varphi \in \Phi _{k}}\, ,\, \left\{ E^{\varphi }\right\} _{\varphi \in \Phi _{k}}$
are pairwise disjoint. For each $\varphi \in \Phi _{k}$ we have by
the above that either (a) $|E^{\varphi }\cap D^{\varphi }|\geq \frac{1}{4C}|E^{\varphi }|$
or (b) $|D^{\varphi }|\geq \frac{1}{4C}2^{(\varepsilon /2)|F_{k}|}|E^{\varphi }|$.
Let $\Phi _{k}^{+}\subseteq \Phi _{k}$ be the set of $\varphi $'s
for which the first alternative holds, and $\Phi _{k}^{-}$ its complement
in $\Phi _{k}$. There are again two alternatives: 
\begin{enumerate}
\item $|\cup _{\varphi \in \Phi _{k}^{+}}E^{\varphi }|\geq \frac{1}{2}|\cup _{\varphi \in \Phi _{k}}E^{\varphi }|$.
In this case\[
|\bigcup _{\varphi \in \Phi _{k}}(E^{\varphi }\cap D^{\varphi })|\geq |\bigcup _{\varphi \in \Phi _{k}^{+}}(E^{\varphi }\cap D^{\varphi })|=\]
 \[
=\sum _{\varphi \in \Phi _{k}^{+}}\left|E^{\varphi }\cap D^{\varphi }\right|\geq \frac{1}{4C}\sum _{\varphi \in \Phi _{k}^{+}}|E^{\varphi }|=\]
\[
=\frac{1}{4C}|\bigcup _{\varphi \in \Phi _{k}^{+}}E^{\varphi }|\geq \frac{1}{8C}|\bigcup _{\varphi \in \Phi _{k}}E^{\varphi }|.\]

\item $|\cup _{\varphi \in \Phi _{k}^{-}}E^{\varphi }|>\frac{1}{2}|\cup _{\varphi \in \Phi _{k}}E^{\varphi }|$.
In this case\[
|\bigcup _{\varphi \in \Phi _{k}}D^{\varphi }|\geq |\bigcup _{\varphi \in \Phi _{k}^{-}}D^{\varphi }|=\]
\[
=\sum _{\varphi \in \Phi _{k}^{-}}\left|D^{\varphi }\right|\geq \frac{1}{4C}2^{(\varepsilon /2)|F_{k}|}\sum _{\varphi \in \Phi _{k}^{-}}|E^{\varphi }|=\]
 \[
\geq \frac{1}{4C}2^{(\varepsilon /2)|F_{k}|}|\bigcup _{\varphi \in \Phi _{k}^{+}}E^{\varphi }|\geq \frac{1}{8C}2^{(\varepsilon /2)|F_{k}|}|\bigcup _{\varphi \in \Phi _{k}}E^{\varphi }|.\]

\end{enumerate}
\noindent \textbf{Step 3:} We now let $k$ vary between $L$ and $M$.
Write $E_{k}=\cup _{\varphi \in \Phi _{k}}E^{\varphi }$ and $D_{k}=\cup _{\varphi \in \Phi _{k}}D^{\varphi }$.
Note that\[
E_{k}=\left\{ f\in F_{N}\, :\, fx\in B\textrm{ and }k(f)=k\right\} \]
 \[
D_{k}=\left\{ f\in F_{N}\, :\, fx\in [\Phi _{k}]\textrm{ and }T_{k}^{(F,Y)}(fx)<2^{(h-\varepsilon /2)|F_{k}|}\right\} \]
 so $\left\{ E_{k}\right\} _{L\leq k\leq M}$ are pairwise disjoint,
but $\left\{ D_{k}\right\} _{L\leq k\leq M}$ need not be. We saw
in step (2) that for each $L\leq k\leq M$, either (a) $|E_{k}\cap D_{k}|\geq \frac{1}{8C}|E_{k}|$
or (b) $|D_{k}|\geq \frac{1}{8C}2^{(\varepsilon /2)|F_{k}|}|E_{k}|$.
Write \[
J=\left\{ k\, :\, L\leq k\leq M\textrm{ and }|D_{k}|\geq \frac{1}{8C}2^{(\varepsilon /2)|F_{k}|}|E_{k}|\right\} .\]
 Let $k\in J$. If $|E_{k}|\geq 16C\cdot 2^{-(\varepsilon /2)|F_{k}|}|F_{N}|$
we would have \[
|D_{k}|\geq \frac{1}{8C}2^{(\varepsilon /2)|F_{k}|}|E_{k}|\geq 2|F_{N}|.\]
 However, $D_{k}\subseteq \left(\cup _{L\leq k\leq M}W_{k}\right)|F_{N}|$
and the latter set, if $N$ is large enough, has size strictly less
than $2|F_{N}|$, so the inequality $|D_{k}|\geq 2|F_{N}|$ is impossible.
Thus for $N$ large enough, if $k\in J$ it must be that $|E_{k}|\leq 16C\cdot 2^{-(\varepsilon /2)|F_{k}|}|F_{N}|$,
so\[
\frac{1}{|F_{N}|}|\bigcup _{k\in J}E_{k}|=\frac{1}{|F_{N}|}\sum _{k\in J}|E_{k}|\leq \frac{16C}{|F_{n}|}\sum _{k=L}^{M}2^{-(\varepsilon /2)|F_{k}|}\leq \frac{p}{2C|F_{N}|}\]
 for $L$ large enough. Therefore,\[
\frac{1}{|F_{N}|}|\bigcup _{L\leq k\leq M}D_{k}|\geq \frac{1}{|F_{N}|}|\bigcup _{k\notin J}E_{k}\cap D_{k}|=\frac{1}{|F_{N}|}\sum _{k\notin J}|E_{k}\cap D_{k}|\geq \]
 \[
\geq \frac{1}{8C|F_{N}|}\sum _{k\notin J}|E_{k}|=\frac{1}{8C|F_{N}|}|\bigcup _{k\notin J}E_{k}|=\]
 \[
=\frac{1}{8C|F_{N}|}|\bigcup _{L\leq k\leq M}E_{k}|-\frac{1}{8C|F_{N}|}|\bigcup _{k\in J}E_{k}|\geq \]
 by assumption, $\frac{1}{|F_{N}|}\left|\bigcup _{L\leq k\leq M}E_{k}\right|\geq p$,
so\[
\geq \frac{p}{8C}-\frac{p}{16C}=\frac{p}{16C}.\]
 But this is impossible, since it implies that at least a $\frac{p}{16C}$-fraction
of $f\in F_{N}$ are such that $\inf _{L\leq k\leq M}T_{k}^{(F,Y)}(fx)<2^{(h-\varepsilon /2)|F_{k}|}$.
Since $L$ was arbitrarily large, this contradicts the fact that $T_{*}^{(F,Y)}\geq h$
a.s.. This completes the proof. 
\end{proof}
In order to prove the theorem for $\{W_{n}\}$ a quasi-incompressible
sequence (definition \ref{def:quasi-incompressible}) we need an analogue
of lemma \ref{lemma:incompressibility}. This lemma is used in the
proof in step (1) to justify equations \ref{eq:First-Use-Of-Incompressibility}
and \ref{eq:Second-Use-Of-Incompressibility}. Recalling what occurred
there, we had fixed $k$ and $\varphi \in \Sigma ^{F_{k}}$ and considered
a set $E\subseteq F_{N}$ of points $f$ such that $U_{k,n(f)}^{(F,W)}(fx)<2^{-(h-\varepsilon )|F_{k}|}$.
Note that the right-hand side of this inequality depends only on $k$.
Furthermore, one of our assumptions was that $|W_{n(f)}|>2^{(h-\varepsilon /2)|F_{k}|}$
for $f\in E$. It is not hard to see that under these conditions a
version of lemma \ref{lemma:incompressibility} is valid, assuming
that the $W_{n-1}$-boundary of $W_{n}$ is small, in a manner depending
on $k$, for all $n$ such that $|W_{n}|>2^{(h-\varepsilon /2)|F_{k}|}$,
and this is just what quasi-incompressibility means. We omit the details.

\section{\label{sec:open-questions}Open questions}

In conclusion we would like mention several problems which remain
unresolved.

The method of proof of the fact that $T_{*}\geq h$ made it necessary
that in the definition of $R_{k}$ we require that a repetition of
the pattern $x(F_{k})$ be {}``counted'' only if it is disjoint
from the original pattern. This definition affects also the definition
of $U_{k,n}$ and the results about $U_{*}$. In the case of $G=\mathbb{Z}^{d}$
studied in \cite{OW02} no such restriction was necessary.

\begin{question}

If we define \[ R^{(F,W)}_{k}(x)=\inf \{n\,:\,\textrm{there exists }1_{G}\neq f\in W_{n} \textrm{ s.t }fx\in x(F_{k}) \} \] and \(T_{k}=|W_{R_{k}}|\), when is it true that \(\frac{1}{|F_{k}|}\log T_{k} \geq h\) ?

\end{question}

For the upper bounds $T^{*}\leq h$ and $U^{*}\leq h$ to hold, we
saw that it is necessary to require some special properties of the
window sets, eg that they be (quasi) filling. We have seen that such
sequences exist in many cases, but we have no information about the
case of non-locally-finite torsion groups:

\begin{question}
Do all groups possess (quasi) filling/incompressible sequences? What classes of groups have (quasi) filling/incompressible F\o{}lner sequences?
\end{question}

Another issue is, under what conditions, milder than those given,
does the bound $U_{*}\geq h$ hold? For one thing, the dependence
on an interpolation sequence would seem to be an artifact of the proof,
and perhaps can be removed. It is also not impossible that the requirement
that $\left\{ W_{n}\right\} $ be incompressible might be weakened
to being filling, or perhaps even completely done away with; the fact
that $T_{*}\geq h$ for very general window sequences, and that we
have not yet found a counterexample of the type in section \ref{subsec:bad-window-sequences}
for the lower bound $U_{*}\geq h$, gives some hope that this might
be true.

\begin{question}
Can weaker conditions on the window sequence be found which ensure the bound \(U_{*}\geq h\)?
\end{question}

\bibliographystyle{plain}
\bibliography{bib}

\end{document}